\documentclass[a4paper,12pt]{article}
\usepackage[hiresbb]{graphicx}
\usepackage{amsmath,amsthm,amssymb}
\usepackage{natbib,color,bm}
\numberwithin{equation}{section}
\usepackage{setspace} 
\usepackage{diagbox}
\usepackage[top=30truemm,bottom=30truemm,left=25truemm,right=25truemm]{geometry}
\usepackage{comment}
\usepackage{algorithmicx}%
\usepackage{algpseudocode}%
\usepackage{listings}%
\usepackage{anyfontsize}%
\usepackage{diagbox}
\usepackage{comment}
\usepackage{booktabs}
\usepackage{authblk}

\usepackage{hyperref}
\hypersetup{%
  colorlinks=true,
  linkcolor=blue,  
  citecolor=black,   
  urlcolor=blue     
}

\makeatletter
\AtBeginDocument{%
  \@ifundefined{NAT@hyper@}{}{%
    \let\NAT@orig@hyper@\NAT@hyper@
    \renewcommand*{\NAT@hyper@}[1]{%
      \begingroup
        \let\NAT@orig@date\NAT@date
        \def\NAT@date{\textcolor{blue}{\NAT@orig@date}}%
        \NAT@orig@hyper@{#1}%
      \endgroup
    }%
  }%
}
\makeatother
 
\theoremstyle{plain} 
\newtheorem{thm}{Theorem}[section]%
\newtheorem{prop}[thm]{Proposition}%
\newtheorem{lem}[thm]{Lemma}%
\newtheorem{cor}[thm]{Corollary}%

\theoremstyle{definition} 
\newtheorem{defn}[thm]{Definition}  
\newtheorem{exam}{Example}
\newtheorem{rem}{Remark}%

\allowdisplaybreaks

\numberwithin{equation}{section}

\def\MC{ {\mathcal{C}} }
\def\MS{ {\mathcal{S}} }

\def\tr{ {\mathrm{tr}}}
\def\ve{ \varepsilon }
\def\rk{ {\mathrm{rank}} }
\def\E{ {\mathrm{E}}}
\def\V{ {\mathrm{Cov}}}
\def\fa{ {\rm for \ all} }

\def\say{ {\rm say} }
\def\an{ {\rm and} }
\def\fs{ {\rm for \ some} }

\def\diag{ {\rm diag} }
\def\MSEM{ {\rm MSEM} }
\def\BRM{ {\rm BRM} }
\def\BRE{ {\rm BRE} }

\title{\Large\textbf{Bayes linear estimator in the general linear model}}
\author[1,2]{Hirai Mukasa}

\affil[1]{Graduate School of Mathematics, Kyushu University}
\affil[2]{National Institute of Technology, Kagoshima College}

\date{}

\begin{document}
\maketitle

\begin{abstract}
The Bayes linear estimator is obtained by minimizing the Bayes risk matrix under squared loss among all linear estimators.
In this paper, we study the statistical properties and equivalence problems of Bayes linear estimators in the general linear model.
First, we examine linear sufficiency and linear completeness of Bayes linear estimators.
Second, we derive necessary and sufficient conditions under which two Bayes linear estimators coincide. 
These conditions clarify when a Bayes linear estimator based on a simpler covariance structure retains Bayes-risk optimality under the original covariance structure.
Moreover, several examples, including Rao’s mixed-effects model and the spatial error model, show that our results can simplify the estimation procedure.
Finally, we establish equivalent conditions for the equality of residual sums of squares when Bayes linear estimators are considered.
\end{abstract}

\maketitle

\section{Introduction}\label{sec:1}	
Consider the general linear model
\begin{align}\label{glm}
 \boldsymbol{y} = \boldsymbol{X} \boldsymbol{\beta} +\boldsymbol{u}, \quad \E[\boldsymbol{u}] = \boldsymbol{0}_n, \quad \V(\boldsymbol{u}) = \sigma^2 \boldsymbol{\Omega},  
\end{align}
where $\boldsymbol{y}\in \mathbb{R}^n$ is an observable random vector, $\boldsymbol{X} \in \mathbb{R}^{n \times k}$ is a non-stochastic matrix $(n>k)$ satisfying $\rk(\boldsymbol{X}) =k$, $\boldsymbol{\beta} \in \mathbb{R}^k$ is an unknown vector, $\sigma^2 > 0$ is an unknown  constant, and $\boldsymbol{\Omega} \in \MS^{+}(n)$ is a covariance matrix.
Here, $\MS^+(m)$ and $\MS^N(m)$ denote the sets of $m \times m$ positive definite and positive semidefinite matrices, respectively.

The ordinary least squares estimator (OLSE) 
\begin{align*}
\hat{\boldsymbol{\beta}}_{OLS} = (\boldsymbol{X}^\top \boldsymbol{X})^{-1}\boldsymbol{X}^\top \boldsymbol{y}
\end{align*}
and the generalized least squares estimator (GLSE) 
\begin{align*}
\hat{\boldsymbol{\beta}}_{GLS} = (\boldsymbol{X}^\top \boldsymbol{\Omega}^{-1}\boldsymbol{X})^{-1}\boldsymbol{X}^\top \boldsymbol{\Omega}^{-1}\boldsymbol{y}
\end{align*}
have long been used as fundamental tools for estimating $\boldsymbol{\beta}$ in \eqref{glm}.
Let $\boldsymbol{L}\boldsymbol{y}$ be a linear estimator of $\boldsymbol{\beta}$.
When $\boldsymbol{\Omega}$ is known, GLSE is the best linear unbiased estimator (BLUE) for $\boldsymbol{\beta}$, that is, $\V(\hat{\boldsymbol{\beta}}_{GLS}) \leq_\mathrm{L} \V(\boldsymbol{L}\boldsymbol{y})$
for any $\boldsymbol{L} \in \mathbb{R}^{k \times n}$ such that $\boldsymbol{L}\boldsymbol{X} = \boldsymbol{I}_k$, where $\boldsymbol{I}_m \in \MS^+(m)$ denotes the $m \times m$ identity matrix.
Here, for symmetric matrices $\boldsymbol{A},\boldsymbol{B}\in \mathbb{R}^{m \times m}$, $\boldsymbol{A} \leq_\mathrm{L} \boldsymbol{B}$ means $\boldsymbol{B} - \boldsymbol{A}\in \mathcal{S}^N(m)$.
However, when there is strong multicollinearity among explanatory variables, the covariance matrix of unbiased estimators becomes large, and such estimators are not necessarily desirable from the viewpoint of mean squared error. Therefore, biased linear estimators such as the ordinary ridge estimator \citep{RefHK70a, RefHK70b} and shrinkage estimators \citep{RefMW73} have been proposed.
The term ``biased estimator" refers to an estimator that is not unbiased for $\boldsymbol{\beta}$.

In this paper, we focus on the linear Bayesian method proposed by \cite{RefH69}, which is a kind of biased linear estimation.
In recent years, linear Bayesian methods have been applied to linear regression models with equality and inequality constraints, as well as to measurement error models \citep{RefZW21,RefJ22,RefW24}.
Let $\pi$ be a prior distribution on $(\boldsymbol{\beta}, \sigma^2)$ satisfying $\E_\pi[\boldsymbol{\beta}\boldsymbol{\beta}^\top] = \boldsymbol{W}, \  \E_\pi[\sigma^2] = \gamma$.
Let us denote $\boldsymbol{K}=\gamma^{-1}\boldsymbol{W}$.
Then the Bayes linear estimator (BLE) 
\begin{align}\label{ble}
\hat{\boldsymbol{\beta}}_{BL}(\boldsymbol{\Omega}, \boldsymbol{K}) = \boldsymbol{K}\boldsymbol{X}^\top(\boldsymbol{\Omega}+\boldsymbol{X}\boldsymbol{K}\boldsymbol{X}^\top)^{-1}\boldsymbol{y}
\end{align}
is derived as a linear estimator that minimizes the Bayes risk matrix in \eqref{glm}
\begin{align}\label{BRML}
\BRM(\boldsymbol{L}\boldsymbol{y}) = \E_\pi[\E[(\boldsymbol{L}\boldsymbol{y} - \boldsymbol{\beta})(\boldsymbol{L}\boldsymbol{y} - \boldsymbol{\beta})^\top \ | \ \boldsymbol{\beta}, \sigma^2]]
\end{align}
with respect to the L\"{o}wner partial ordering, that is, $\BRM(\hat{\boldsymbol{\beta}}_{BL}(\boldsymbol{\Omega}, \boldsymbol{K})) \leq_\mathrm{L} \BRM(\boldsymbol{L}\boldsymbol{y})$ for any $\boldsymbol{L} \in \mathbb{R}^{k \times n}$.  
Thus, the Bayes risk matrix provides a criterion for evaluating linear estimators by incorporating both covariance matrix and bias vector. 
This is one of the main motivations for studying BLE under \eqref{glm}.

Bayes estimators are obtained as posterior means under squared loss and therefore involve integral calculations, which often makes it difficult to obtain closed-form expressions. 
In contrast, BLE depends only on the second-order moments of the prior distribution and can be expressed in explicit form.
Several important properties of BLE have been studied in the literature.
\cite{RefR76} showed that the extended Bayes class, which consists of BLEs and their limit points, coincides with the class of linearly admissible estimators.
He also noted that when $\boldsymbol{K} \in \MS^+(k)$, BLE can be expressed as a general ridge estimator (GRE) introduced in \cite{RefR75}
and discussed its relation to minimaxity. 
\cite{RefL78} pointed out that the extended Bayes class includes not only unbiased estimators but also general ridge, fractional rank \citep{RefM70}, and certain restricted least squares estimators \citep{RefR07}. 
\cite{RefH96} investigated choices of weight matrices for which the BLE is minimax with respect to generalized quadratic risk functions.
\cite{RefG03} also discussed quasi minimaxity of BLE in certain situations and showed that BLE coincides with GLSE under a noninformative prior distribution. 
For linear admissibility and minimaxity, see, for example, \cite{RefB75,RefBM88,RefB89,RefGM04,RefAS00,RefAS11}. A summary of properties of BLE is also given in \cite{RefM26}.

The aim of this paper is to study statistical properties and equivalence problems of BLEs in \eqref{glm}.
In Section~\ref{sec:2}, we derive BLE and introduce some preliminary results used in the subsequent proofs.
In Section~\ref{sec:3}, we introduce definitions of linear sufficiency and linear completeness and examine them for BLE.
These properties characterize whether a linear estimator retains sufficient and
nonredundant information to construct the BLUE for $\boldsymbol{X}\boldsymbol{\beta}$.
Related work has developed these concepts in several classes of linear models, including continuous-time Gauss--Markov models \citep{RefIP03,RefIP04} and the partitioned linear model \citep{RefIP06}. 
In addition, \cite{RefK17} summarized the equivalence of BLUEs under the general linear model and its linearly transformed model from the viewpoint of linear sufficiency, and \cite{RefH20} examined the relationship between linear sufficiency and best linear unbiased predictor in the linear mixed model.
These studies mainly concern linear transformations arising in various linear models, but it is not immediately clear whether a certain linear estimator possesses linear sufficiency or linear completeness. 
Therefore, we examine these properties for BLE and clarify how they depend on $\boldsymbol{K}$.

For an observation $\boldsymbol{y}$, we consider the following two models with different prior distributions $\pi_1$ and $\pi_2$:
\begin{align}\label{model}
\left\{ \,
 \begin{aligned}
& \textup{Model 1:}  \ \boldsymbol{y} \ | \ \boldsymbol{\beta}, \sigma^2 \sim (\boldsymbol{X}\boldsymbol{\beta}, \sigma^2 \boldsymbol{\Omega}), \ \E_{\pi_1}[\boldsymbol{\beta}\boldsymbol{\beta}^\top] = \boldsymbol{W}_1, \ \E_{\pi_1}[\sigma^2] = \gamma_1,\\
&  \textup{Model 2:}  \ \boldsymbol{y} \ | \ \boldsymbol{\beta}, \sigma^2 \sim (\boldsymbol{X}\boldsymbol{\beta}, \sigma^2 \boldsymbol{I}_n), \ \E_{\pi_2}[\boldsymbol{\beta}\boldsymbol{\beta}^\top] = \boldsymbol{W}_2, \ \E_{\pi_2}[\sigma^2] = \gamma_2,
\end{aligned}
    \right.
\end{align}
where $\boldsymbol{y}\mid\boldsymbol{\beta},\sigma^2\sim(\boldsymbol{\mu},\boldsymbol{\Sigma})$
means that the conditional mean and covariance matrix of $\boldsymbol{y}$ given $\boldsymbol{\beta}$ and $\sigma^2$ are $\boldsymbol{\mu}$ and $\boldsymbol{\Sigma}$, respectively.
The BLEs under \eqref{model} are obtained as
\begin{align*}
&\hat{\boldsymbol{\beta}}_{BL}(\boldsymbol{\Omega}, \boldsymbol{K}_1) = \boldsymbol{K}_1\boldsymbol{X}^\top (\boldsymbol{\Omega} + \boldsymbol{X}\boldsymbol{K}_1\boldsymbol{X}^\top)^{-1}\boldsymbol{y},\\
&\hat{\boldsymbol{\beta}}_{BL}(\boldsymbol{I}_n, \boldsymbol{K}_2) = \boldsymbol{K}_2\boldsymbol{X}^\top (\boldsymbol{I}_n + \boldsymbol{X}\boldsymbol{K}_2\boldsymbol{X}^\top)^{-1}\boldsymbol{y},
\end{align*}
respectively, where $\boldsymbol{K}_i=\gamma_i^{-1}\boldsymbol{W}_i$ is fixed for $i =1,2$.
In Section~\ref{sec:4.1}, we derive necessary and sufficient conditions under which the identical equality
\begin{align}\label{B2E}
\hat{\boldsymbol{\beta}}_{BL}(\boldsymbol{\Omega}, \boldsymbol{K}_1) \equiv \hat{\boldsymbol{\beta}}_{BL}(\boldsymbol{I}_n, \boldsymbol{K}_2) 
\end{align}
holds.
Here, $\equiv$ means that the equality holds for all $\boldsymbol{y}\in\mathbb{R}^n$.
This problem is worth considering from both theoretical and practical viewpoints.
If \eqref{B2E} holds, then BLE derived under Model 2 also minimizes the Bayes risk matrix under Model 1. Thus, the optimality of $\hat{\boldsymbol{\beta}}_{BL}(\boldsymbol{I}_n,\boldsymbol{K}_2)$ is retained even when the true covariance structure is $\boldsymbol{\Omega}$ rather than $\boldsymbol{I}_n$.
In related work, \cite{RefMM73} and \cite{RefH12} compared BLUEs across linear models with different covariance structures, and \cite{RefHIMP23,RefHIMP25} investigated classes of covariance structures under which BLUEs are retained.
Thus, the present study extends this line of research to best linear estimators obtained by minimizing the Bayes risk matrix.
Second, the equivalence can simplify estimation procedures.
When $\boldsymbol{\Omega}$ is partially unknown, we often estimate $\boldsymbol{\beta}$ by using a two-step estimation procedure in which we first construct an estimator $\hat{\boldsymbol{\Omega}}$ and then compute $\hat{\boldsymbol{\beta}}_{BL}(\hat{\boldsymbol{\Omega}},\boldsymbol{K}_1)$.
However, if the equality \eqref{B2E} is satisfied, then it is sufficient to use $\hat{\boldsymbol{\beta}}_{BL}(\boldsymbol{I}_n, \boldsymbol{K}_2)$, and the two-step estimation is unnecessary.
In Section~\ref{sec:4.2}, we present practical applications and verify the results through numerical experiments.

Now, related problems have been studied extensively in the literature.
In particular, conditions under which OLSE coincides with GLSE have been investigated by many authors; see, for example, \cite{RefA48, RefZ67, RefK68, RefPS89}.
For a recent overview,  see \cite{RefM21}.
Among others, \cite{RefR67} and \cite{RefG70} showed that
the equality $\hat{\boldsymbol{\beta}}_{GLS} \equiv \hat{\boldsymbol{\beta}}_{OLS}$ holds if and only if $\boldsymbol{\Omega}$ is of the form
\begin{align}\label{Rcond}
\boldsymbol{\Omega} = \boldsymbol{X} \boldsymbol{\Gamma} \boldsymbol{X}^\top + \boldsymbol{Z} \boldsymbol{\Delta} \boldsymbol{Z}^\top \quad  \fs \ \ \boldsymbol{\Gamma} \in \mathcal{S}^+(k), \  \boldsymbol{\Delta} \in \mathcal{S}^+(n-k),
\end{align}
where $\boldsymbol{Z} \in \mathbb{R}^{n \times (n-k)}$ is any matrix such that $\boldsymbol{X}^\top \boldsymbol{Z} = \boldsymbol{0}$ and $\rk(\boldsymbol{Z}) = n-k$.
The structure \eqref{Rcond} of $\boldsymbol{\Omega}$ is called
Rao's covariance structure.
In addition, \cite{RefTW06} studied equality and proportionality among OLSE, the weighted least squares estimator (WLSE), and GLSE, while \cite{RefLP11} applied conditions for the equivalence between WLSE and GLSE to local polynomial regression and kernel smoothing.
\cite{RefTK20} and \cite{RefMT25} derived equivalent conditions under which two GREs coincide.
Although equivalence problems have been studied for unbiased linear estimators and GREs, no existing studies have explored the equivalence of BLEs. 
Thus, this paper fills this gap.

In Section~\ref{sec:5},  instead of looking for $(\boldsymbol{\Omega},\boldsymbol{X})$ such that two BLEs coincide for all $\boldsymbol{y} \in \mathbb{R}^n$, we characterize the set of $\boldsymbol{y}$ satisfying
\begin{align}\label{B2Ey}
\hat{\boldsymbol{\beta}}_{BL}(\boldsymbol{\Omega}, \boldsymbol{K}_1) = \hat{\boldsymbol{\beta}}_{BL} (\boldsymbol{I}_n, \boldsymbol{K}_2) 
\end{align}
for given $(\boldsymbol{\Omega}, \boldsymbol{X})$.
In previous studies, \cite{RefK80b} proved that $\hat{\boldsymbol{\beta}}_{GLS} = \hat{\boldsymbol{\beta}}_{OLS}$ if and only if
\begin{align*}
\boldsymbol{y} \in \MC(\boldsymbol{X}) \oplus \left[\MC(\boldsymbol{X})^\bot\ \cap \ \MC(\boldsymbol{\Omega}\boldsymbol{Z})\right],
\end{align*}
where $\MC(\boldsymbol{M})$ denotes the column space of a matrix $\boldsymbol{M}$, and for a subspace $W\subseteq \mathbb{R}^n$, $W^{\perp}$ denotes its orthogonal complement, that is, $\mathbb{R}^n = W \oplus W^{\perp}$  with a direct sum $\oplus$.
For a similar study, see \cite{RefGT97} and \cite{RefG01b}.
Since BLE is biased for $\boldsymbol{\beta}$, the condition for the equality \eqref{B2Ey} cannot be expressed in terms of a direct sum decomposition.
However, our result shows that the set of $\boldsymbol{y}$ for the equality \eqref{B2Ey} can be identified in this form by imposing \eqref{Rcond}.

If we consider a quadratic estimator of $\sigma^2$, we use the residual sums of squares of $\hat{\boldsymbol{\beta}}_{BL}(\boldsymbol{\Phi},\boldsymbol{K})$
\begin{align*}
RSS_{BL}(\boldsymbol{\Phi},\boldsymbol{K}) = \left\| \boldsymbol{\Phi}^{-1/2} \left( \boldsymbol{y} - \boldsymbol{X} \hat{\boldsymbol{\beta}}_{BL}(\boldsymbol{\Phi}, \boldsymbol{K}) \right) \right\|^2, \quad \boldsymbol{\Phi} \in \{\boldsymbol{\Omega}, \boldsymbol{I}_n\}, \
\boldsymbol{K} \in \MS^N(k),
\end{align*}
where $\|\cdot\|$ stands for the Euclidean norm.
In Section~\ref{sec:6}, we derive necessary and sufficient conditions under which the identical equality
\begin{align} \label{BRSS1}
 RSS_{BL}(\boldsymbol{\Omega},\boldsymbol{K}_1) \equiv RSS_{BL}(\boldsymbol{I}_n,\boldsymbol{K}_2)
\end{align}
holds.
In particular, \cite{RefK68} and \cite{RefK80} studied the equivalence of residual sums of squares associated with OLSE and GLSE. More recently, \cite{RefMT25} derived equivalent conditions under which two residual sums of squares of GREs are equivalent.
Note that the necessary and sufficient condition for the equality \eqref{B2E} is broader than the condition that $\boldsymbol{\Omega} = \boldsymbol{I}_n$ and $\boldsymbol{K}_1 = \boldsymbol{K}_2$.
We therefore ask what additional conditions, together with the equality \eqref{B2E}, guarantee that this trivial condition must hold.
To answer this question, we consider the case where the two equalities \eqref{B2E} and \eqref{BRSS1} simultaneously hold.

Finally, some concluding remarks are presented in Section~\ref{sec:7}.

\section{Preliminaries}\label{sec:2}
In this section, we derive BLE in \eqref{glm}  and introduce some preliminary results  used in subsequent sections. 

The mean squared error matrix of a linear estimator $\boldsymbol{L}\boldsymbol{y}$ can be written as
\begin{align*}
\MSEM(\boldsymbol{L}\boldsymbol{y}) &= \E[(\boldsymbol{L}\boldsymbol{y} - \boldsymbol{\beta})(\boldsymbol{L}\boldsymbol{y} - \boldsymbol{\beta})^\top \ | \ \boldsymbol{\beta}, \sigma^2]\\
&= \V(\boldsymbol{L}\boldsymbol{y}) + (\E[\boldsymbol{L}\boldsymbol{y}] - \boldsymbol{\beta})(\E[\boldsymbol{L}\boldsymbol{y}] - \boldsymbol{\beta})^\top \\
&= \sigma^2 \boldsymbol{L}\boldsymbol{\Omega}\boldsymbol{L}^\top + (\boldsymbol{L}\boldsymbol{X} - \boldsymbol{I}_k) \boldsymbol{\beta}\boldsymbol{\beta}^\top (\boldsymbol{L}\boldsymbol{X} - \boldsymbol{I}_k)^\top.
\end{align*}
Since $\boldsymbol{\beta}$ and $\sigma^2$ are unknown, this expression cannot be used directly to evaluate the performance of $\boldsymbol{L}\boldsymbol{y}$.
To obtain a workable criterion, we adopt a Bayesian approach by placing a prior distribution $\pi$ on $(\boldsymbol{\beta}, \sigma^2)$ that satisfies $\E_{\pi}[\boldsymbol{\beta}\boldsymbol{\beta}^\top] = \boldsymbol{W}$ and $\E_{\pi}[\sigma^2] = \gamma$.
Then \eqref{BRML} can be expressed as
\begin{align}\label{BRM1}
\BRM(\boldsymbol{L}\boldsymbol{y}) &= \E_\pi[\E[(\boldsymbol{L}\boldsymbol{y} - \boldsymbol{\beta})(\boldsymbol{L}\boldsymbol{y} - \boldsymbol{\beta})^\top \ | \ \boldsymbol{\beta}, \sigma^2]] \notag\\
& = \gamma \boldsymbol{L}\boldsymbol{\Omega}\boldsymbol{L}^\top + (\boldsymbol{L}\boldsymbol{X} - \boldsymbol{I}_k) \boldsymbol{W} (\boldsymbol{L}\boldsymbol{X} - \boldsymbol{I}_k)^\top \notag\\
& = \gamma \ \{\boldsymbol{L}(\boldsymbol{\Omega} + \boldsymbol{X}\boldsymbol{K}\boldsymbol{X}^\top)\boldsymbol{L}^\top - \boldsymbol{L}\boldsymbol{X}\boldsymbol{K}-\boldsymbol{K}\boldsymbol{X}^\top\boldsymbol{L}^\top + \boldsymbol{K}\},
\end{align}
where $\boldsymbol{K} = \gamma^{-1}\boldsymbol{W}$.

Proposition~\ref{pr21} provides a general expression for the Bayes risk matrix and an explicit form of BLE that minimizes the Bayes risk matrix.
The former will be used in the proof of Proposition~\ref{pr45}.
\begin{prop}\label{pr21}
The Bayes linear estimator that minimizes \eqref{BRML} in \eqref{glm} with respect to the L\"{o}wner partial ordering is given in the form of \eqref{ble}.
Moreover, the Bayes risk matrix for $\boldsymbol{L}\boldsymbol{y}$ can be expressed as
\begin{align}\label{BRM2}
\BRM(\boldsymbol{L}\boldsymbol{y}) =  \gamma (\boldsymbol{L} - \boldsymbol{L}_*)(\boldsymbol{\Omega} + \boldsymbol{X}\boldsymbol{K}\boldsymbol{X}^\top)(\boldsymbol{L} - \boldsymbol{L}_*)^\top + \BRM(\boldsymbol{L}_*\boldsymbol{y}),
\end{align}
where $\boldsymbol{L}_* = \boldsymbol{K}\boldsymbol{X}^\top(\boldsymbol{\Omega}+\boldsymbol{X}\boldsymbol{K}\boldsymbol{X}^\top)^{-1}$.
\end{prop}
\begin{proof}
To simplify the notation, let $\boldsymbol{T} = \boldsymbol{\Omega} + \boldsymbol{X}\boldsymbol{K}\boldsymbol{X}^\top$.
The first term on \eqref{BRM2} can be rewritten as
\begin{align*}
(\boldsymbol{L} - \boldsymbol{L}_*)\boldsymbol{T}(\boldsymbol{L} - \boldsymbol{L}_*)^\top
& = \boldsymbol{L}\boldsymbol{T}\boldsymbol{L}^\top - \boldsymbol{L}\boldsymbol{T}\boldsymbol{L}_*^\top - \boldsymbol{L}_* \boldsymbol{T}\boldsymbol{L}^\top + \boldsymbol{L}_*\boldsymbol{T}\boldsymbol{L}_*^\top\\
& = \boldsymbol{L}\boldsymbol{T}\boldsymbol{L}^\top -\boldsymbol{L}\boldsymbol{X}\boldsymbol{K} - \boldsymbol{K}\boldsymbol{X}^\top \boldsymbol{L}^\top + \boldsymbol{K}\boldsymbol{X}^\top\boldsymbol{T}^{-1}\boldsymbol{X}\boldsymbol{K}.
\end{align*}   
From \eqref{BRM1}, the second term on \eqref{BRM2} can be calculated as
\begin{align*}
\gamma^{-1}\BRM(\boldsymbol{L}_*\boldsymbol{y}) &= \boldsymbol{L}_*\boldsymbol{T}\boldsymbol{L}_*^\top - \boldsymbol{L}_*\boldsymbol{X}\boldsymbol{K}-\boldsymbol{K}\boldsymbol{X}^\top\boldsymbol{L}_*^\top + \boldsymbol{K}\\
& = \boldsymbol{K}\boldsymbol{X}^\top\boldsymbol{T}^{-1}\boldsymbol{X}\boldsymbol{K} - \boldsymbol{K}\boldsymbol{X}^\top\boldsymbol{T}^{-1}\boldsymbol{X}\boldsymbol{K} - \boldsymbol{K}\boldsymbol{X}^\top\boldsymbol{T}^{-1}\boldsymbol{X}\boldsymbol{K} + \boldsymbol{K}\\
& = \boldsymbol{K} - \boldsymbol{K}\boldsymbol{X}^\top\boldsymbol{T}^{-1}\boldsymbol{X}\boldsymbol{K}.
\end{align*}
Therefore, \eqref{BRM1} coincides with \eqref{BRM2}.

Since $\boldsymbol{T} \in \MS^+(n)$, the first term on \eqref{BRM2} is nonnegative definite, and equals zero if and only if $\boldsymbol{L} = \boldsymbol{L}_*$.
Thus, the Bayes linear estimator is given in the form of \eqref{ble}. 
This completes the proof.
\end{proof}

Next, the following lemma provides an alternative representation of BLE.
For a similar proof, see Corollary 3.6 of \cite{RefG03}.
\begin{lem}\label{lem2}
$\hat{\boldsymbol{\beta}}_{BL}(\boldsymbol{\Omega},\boldsymbol{K})$ can be rewritten as
\begin{align}\label{ble2}
\hat{\boldsymbol{\beta}}_{BL}(\boldsymbol{\Omega},\boldsymbol{K}) = \boldsymbol{K}\left[(\boldsymbol{X}^\top \boldsymbol{\Omega}^{-1}\boldsymbol{X})^{-1}+\boldsymbol{K}\right]^{-1}(\boldsymbol{X}^\top \boldsymbol{\Omega}^{-1}\boldsymbol{X})^{-1}\boldsymbol{X}^\top \boldsymbol{\Omega}^{-1}\boldsymbol{y}.
\end{align}
\end{lem}

\begin{rem}
Since the eigenvalues of $\boldsymbol{K}[(\boldsymbol{X}^\top \boldsymbol{\Omega}^{-1}\boldsymbol{X})^{-1}+\boldsymbol{K}]^{-1}$ lie in $[0,1)$, we can interpret BLE as an estimator obtained by
shrinking WLSE toward the origin.
\end{rem}
If $\boldsymbol{K}$ is nonsingular, BLE can be rewritten as follows.
For a similar proof, see Corollary 3.7 of \cite{RefG03}.
\begin{cor}\label{cor23}
If $\boldsymbol{K} \in \MS^+(k)$, $\hat{\boldsymbol{\beta}}_{BL}(\boldsymbol{\Omega},\boldsymbol{K})$ can be written in the form of
\begin{align*}
\hat{\boldsymbol{\beta}}_{BL}(\boldsymbol{\Omega},\boldsymbol{K}) = (\boldsymbol{X}^\top \boldsymbol{\Omega}^{-1}\boldsymbol{X} + \boldsymbol{K}^{-1})^{-1}\boldsymbol{X}^\top \boldsymbol{\Omega}^{-1}\boldsymbol{y}.
\end{align*}
\end{cor}
\begin{rem}\label{rem2}
Corollary~\ref{cor23} implies that when $\boldsymbol{K} \in \MS^+(k)$, $\hat{\boldsymbol{\beta}}_{BL}(\boldsymbol{\Omega},\boldsymbol{K}^{-1})$ coincides with GRE
\begin{align*}
\hat{\boldsymbol{\beta}}_{GR}(\boldsymbol{\Omega}, \boldsymbol{K}) 
 = (\boldsymbol{X}^\top \boldsymbol{\Omega}^{-1}\boldsymbol{X} + \boldsymbol{K})^{-1}\boldsymbol{X}^\top \boldsymbol{\Omega}^{-1}\boldsymbol{y}.
\end{align*}
Therefore, several properties of GRE, such as linear sufficiency and linear admissibility, can be transferred to BLE in such cases.
\end{rem}

Finally, we derive the expectation vector and covariance matrix of BLE.
This fact will be used in the proof of Proposition~\ref{pr47}.
\begin{lem}\label{lem24}
The expectation vector and the covariance matrix of $\hat{\boldsymbol{\beta}}_{BL}(\boldsymbol{\Omega}, \boldsymbol{K})$ in \eqref{glm} are given by
\begin{align*}
& \E\left[\hat{\boldsymbol{\beta}}_{BL}(\boldsymbol{\Omega}, \boldsymbol{K})\right] = \boldsymbol{A}\boldsymbol{\beta}, \quad  \V\left(\hat{\boldsymbol{\beta}}_{BL}(\boldsymbol{\Omega}, \boldsymbol{K})\right) = \sigma^2\boldsymbol{A}(\boldsymbol{X}^\top\boldsymbol{\Omega}^{-1}\boldsymbol{X})^{-1}\boldsymbol{A}^\top,
\end{align*}
where $\boldsymbol{A} = \boldsymbol{K}[(\boldsymbol{X}^\top \boldsymbol{\Omega}^{-1}\boldsymbol{X})^{-1}+\boldsymbol{K}]^{-1}$.
\end{lem}
\begin{proof}
This fact follows immediately from Lemma~\ref{lem2}.
This completes the proof.
\end{proof}
\begin{rem}\label{rm3}
(i) The bias of $\hat{\boldsymbol{\beta}}_{BL}(\boldsymbol{\Omega}, \boldsymbol{K})$ is
\begin{align*}
\E\left[\hat{\boldsymbol{\beta}}_{BL}(\boldsymbol{\Omega}, \boldsymbol{K})\right] - \boldsymbol{\beta}
&= \boldsymbol{K}[(\boldsymbol{X}^\top \boldsymbol{\Omega}^{-1}\boldsymbol{X})^{-1}+\boldsymbol{K}]^{-1}\boldsymbol{\beta} - \boldsymbol{\beta}\\
&= -(\boldsymbol{I}_k + \boldsymbol{K}\boldsymbol{X}^\top\boldsymbol{\Omega}^{-1}\boldsymbol{X})^{-1}\boldsymbol{\beta},
\end{align*}
which implies that BLE is biased for $\boldsymbol{\beta} \ (\neq \boldsymbol{0})$ for any finite choice of $\boldsymbol{K} \in \MS^N(k)$.\\
(ii) The expectation vector and the covariance matrix of $\hat{\boldsymbol{\beta}}_{BL}(\boldsymbol{I}_n, \boldsymbol{K})$ in \eqref{glm} can be calculated as
\begin{align*}
\E\left[\hat{\boldsymbol{\beta}}_{BL}(\boldsymbol{I}_n, \boldsymbol{K})\right] = \boldsymbol{B}\boldsymbol{\beta}, \quad \V\left(\hat{\boldsymbol{\beta}}_{BL}(\boldsymbol{I}_n, \boldsymbol{K})\right) = \sigma^2 \boldsymbol{B}(\boldsymbol{X}^\top \boldsymbol{X})^{-1}\boldsymbol{X}^\top \boldsymbol{\Omega}\boldsymbol{X}(\boldsymbol{X}^\top \boldsymbol{X})^{-1}\boldsymbol{B}^\top,
\end{align*}
where $\boldsymbol{B} = \boldsymbol{K}[(\boldsymbol{X}^\top \boldsymbol{X})^{-1} + \boldsymbol{K}]^{-1}$.
\end{rem}

\section{Linear sufficiency and linear completeness}\label{sec:3}
In this section, we introduce the notions of linear sufficiency and linear completeness and examine whether $\hat{\boldsymbol{\beta}}_{BL}(\boldsymbol{\Omega},\boldsymbol{K})$ possesses these properties.

The definitions of linear sufficiency and linear completeness in the sense of \cite{RefD83} are as follows:
\begin{defn}\label{def1} 
\textup{(i)} A linear transformation  $\boldsymbol{Fy}$ is said to be linearly sufficient if there exists a linear transformation $\boldsymbol{L}$
such that $\boldsymbol{LFy}$ is the BLUE of $\boldsymbol{X}\boldsymbol{\beta}$.\\
\textup{(ii)} A linear transformation $\boldsymbol{Fy}$ is said to be linearly complete if for all linear transformations $\boldsymbol{L}$ such that $\E[\boldsymbol{LFy}] = \boldsymbol{0}$ for all $\boldsymbol{\beta} \in \mathbb{R}^k$,
it follows that $\boldsymbol{LFy} = \boldsymbol{0}$ almost surely.
\end{defn}

Note that if $\boldsymbol{Fy}$ is both linearly sufficient and linearly complete, then the BLUE of $\boldsymbol{X}\boldsymbol{\beta}$ is uniquely determined among linear estimators of the form $\boldsymbol{L}\boldsymbol{F}\boldsymbol{y}$.

To prove Theorem~\ref{thm1}, we present some technical lemmas.
First, we introduce a lemma that rewrites Definition~\ref{def1} in terms of column space relationships; see, for example, \cite{RefM96}.
\begin{lem}\label{lem1}
Let $\boldsymbol{T}= \boldsymbol{\Omega} + \boldsymbol{XWX}^\top$, where 
$\boldsymbol{W}$ is an arbitrary symmetric matrix satisfying $\MC(\boldsymbol{T}) = \MC(\boldsymbol{X}) + \MC(\boldsymbol{\Omega})$.
Then the following holds:\\
(i) $\boldsymbol{Fy}$ is linearly sufficient if and only if $\MC(\boldsymbol{X}) \subseteq \MC(\boldsymbol{T}\boldsymbol{F}^\top)$.\\
(ii) $\boldsymbol{Fy}$ is linearly complete if and only if $\MC(\boldsymbol{F\Omega}) \subseteq \MC(\boldsymbol{FX})$.
\end{lem}
\begin{rem}\label{rm1}
(i) Let $\boldsymbol{A}\in \mathbb{R}^{a \times b}$ and $\boldsymbol{B}\in \mathbb{R}^{a \times c}$.
Then $\MC(\boldsymbol{A}) \subseteq \MC(\boldsymbol{B})$ if and only if there exists a matrix $\boldsymbol{G}\in \mathbb{R}^{c \times b}$ such that $\boldsymbol{A} = \boldsymbol{BG}$.
Moreover, for $\boldsymbol{A}, \boldsymbol{B}\in \mathbb{R}^{a \times b}$, $\MC(\boldsymbol{A}) = \MC(\boldsymbol{B})$ if and only if there exists a nonsingular matrix $\boldsymbol{G} \in \mathbb{R}^{b \times b}$ such that $\boldsymbol{A} = \boldsymbol{BG}$.\\
(ii) Definition~\ref{def1} can also be rephrased as follows.\\
(a) $\boldsymbol{Fy}$ is linearly sufficient if there exists a linear transformation $\boldsymbol{L}$ such that
$\boldsymbol{LFX} = \boldsymbol{X}$ and $\boldsymbol{LF\Omega}\boldsymbol{Z} = \boldsymbol{0}$.\\
(b) $\boldsymbol{Fy}$ is linearly complete if $\boldsymbol{LFX} = \boldsymbol{0}$ implies $\boldsymbol{LF\Omega} = \boldsymbol{0}$ for all linear transformations $\boldsymbol{L}$.
\end{rem}

It is well known that $\hat{\boldsymbol{\beta}}_{GR}(\boldsymbol{\Omega},\boldsymbol{K})$ has both linear sufficiency and linear completeness \citep{RefM96}. 
Since BLE can be expressed in the form of GRE in certain cases (see Remark~\ref{rem2}), it is natural to ask whether $\hat{\boldsymbol{\beta}}_{BL}(\boldsymbol{\Omega},\boldsymbol{K})$ also possesses these properties. The following theorem gives an answer to this question and constitutes one of the main results of this paper.

\begin{thm}\label{thm1}
(i) For $\boldsymbol{K}\in \MS^+(k)$, $\hat{\boldsymbol{\beta}}_{BL}(\boldsymbol{\Omega},\boldsymbol{K})$ is  linearly sufficient.
For $\boldsymbol{K} \in \MS^N(k) \setminus{\MS^+(k)}$, $\hat{\boldsymbol{\beta}}_{BL}(\boldsymbol{\Omega},\boldsymbol{K})$  
is not linearly sufficient.\\
(ii) For $\boldsymbol{K} \in \MS^N(k)$, $\hat{\boldsymbol{\beta}}_{BL}(\boldsymbol{\Omega},\boldsymbol{K})$ is linearly complete.
\end{thm}
\begin{proof}
(i) First, we consider the case where $\boldsymbol{K}\in \MS^+(k)$.
Let $\boldsymbol{T}= \boldsymbol{\Omega} + \boldsymbol{XKX}^\top$ and $\boldsymbol{F}=\boldsymbol{KX}^\top \boldsymbol{T}^{-1}$.
For $\hat{\boldsymbol{\beta}}_{BL}(\boldsymbol{\Omega},\boldsymbol{K})=\boldsymbol{Fy}$, $\boldsymbol{TF}^\top$ is calculated as
\begin{align*}
\boldsymbol{TF}^\top = \boldsymbol{TT}^{-1} \boldsymbol{XK} = \boldsymbol{XK}.
\end{align*}
Let us denote $\boldsymbol{G}=\boldsymbol{K}^{-1}$.
Since $\boldsymbol{K}$ is nonsingular, it holds that $\boldsymbol{X}=\boldsymbol{T}\boldsymbol{F}^\top\boldsymbol{G}$.
Hence, $\MC(\boldsymbol{X}) \subseteq \MC(\boldsymbol{TF}^\top)$ is satisfied.
It follows from Lemma~\ref{lem1}-(i) that $\hat{\boldsymbol{\beta}}_{BL}(\boldsymbol{\Omega},\boldsymbol{K})$ is linearly sufficient.

Next, we consider the case where $\boldsymbol{K} \in \MS^N(k) \setminus{\MS^+(k)}$.
Suppose that $\hat{\boldsymbol{\beta}}_{BL}(\boldsymbol{\Omega},\boldsymbol{K})$ is linearly sufficient.
It follows from Lemma~\ref{lem1}-(i) that there exists a matrix $\boldsymbol{G} \in \mathbb{R}^{k \times k}$ such that $\boldsymbol{X} = \boldsymbol{TF}^\top \boldsymbol{G} =\boldsymbol{XKG}$.
Then premultiplying by $(\boldsymbol{X}^\top \boldsymbol{X})^{-1}\boldsymbol{X}^\top$
yields $\boldsymbol{KG} = \boldsymbol{I}_k$.
The right-hand side has rank $k$, whereas  the left-hand side satisfies $\rk(\boldsymbol{KG}) \leq \rk(\boldsymbol{K})<k$.
This is a contradiction.
Hence, $\hat{\boldsymbol{\beta}}_{BL}(\boldsymbol{\Omega},\boldsymbol{K})$ is not linearly sufficient.\\
(ii) Let $\boldsymbol{F} = \boldsymbol{K}[(\boldsymbol{X}^\top \boldsymbol{\Omega}^{-1}\boldsymbol{X})^{-1}+\boldsymbol{K}]^{-1}(\boldsymbol{X}^\top \boldsymbol{\Omega}^{-1} \boldsymbol{X})^{-1} \boldsymbol{X}^\top \boldsymbol{\Omega}^{-1}$.
From Lemma~\ref{lem2}, we can write $\hat{\boldsymbol{\beta}}_{BL}(\boldsymbol{\Omega},\boldsymbol{K}) = \boldsymbol{Fy}$.
$\boldsymbol{F\Omega}$ and $\boldsymbol{FX}$ are calculated as
\begin{align*}
& \boldsymbol{F\Omega} = \boldsymbol{K}\left[(\boldsymbol{X}^\top \boldsymbol{\Omega}^{-1}\boldsymbol{X})^{-1}+\boldsymbol{K}\right]^{-1}(\boldsymbol{X}^\top \boldsymbol{\Omega}^{-1} \boldsymbol{X} )^{-1}\boldsymbol{X}^\top,\\
& \boldsymbol{FX} = \boldsymbol{K}\left[(\boldsymbol{X}^\top \boldsymbol{\Omega}^{-1}\boldsymbol{X})^{-1}+\boldsymbol{K}\right]^{-1}.
\end{align*}
Letting $\boldsymbol{G} = (\boldsymbol{X}^\top \boldsymbol{\Omega}^{-1} \boldsymbol{X} )^{-1}\boldsymbol{X}^\top$, we obtain $\boldsymbol{F\Omega} = \boldsymbol{FXG}$.
Hence, $\MC(\boldsymbol{F\Omega}) \subseteq \MC(\boldsymbol{FX})$ is satisfied.
It follows from Lemma~\ref{lem1}-(ii) that $\hat{\boldsymbol{\beta}}_{BL}(\boldsymbol{\Omega},\boldsymbol{K})$ is linearly complete.
This completes the proof.
\end{proof}

\begin{rem}
(i) When $\boldsymbol{K} \in \MS^N(k) \setminus{\MS^+(k)}$, 
the linear sufficiency and linear completeness of
$\hat{\boldsymbol{\beta}}_{GR}(\boldsymbol{\Omega},\boldsymbol{K})$ and
$\hat{\boldsymbol{\beta}}_{BL}(\boldsymbol{\Omega},\boldsymbol{K})$ are summarized in Table~\ref{table1}.
\begin{table}[t]
\centering
\setlength{\belowcaptionskip}{8pt}
\caption{Properties of estimators when $\boldsymbol{K} \in \MS^N(k) \setminus{\MS^+(k)}$}
\label{table1}
\begin{tabular}{|c|c|c|}
\hline
\diagbox{Properties}{Estimators}
 & $\hat{\boldsymbol{\beta}}_{GR}(\boldsymbol{\Omega},\boldsymbol{K})$
 & $\hat{\boldsymbol{\beta}}_{BL}(\boldsymbol{\Omega},\boldsymbol{K})$ \\
\hline
Linear sufficiency & $\bigcirc$ & $\times$ \\
\hline
Linear completeness & $\bigcirc$ & $\bigcirc$ \\
\hline
\end{tabular}%
\end{table}
We find that only $\hat{\boldsymbol{\beta}}_{GR}(\boldsymbol{\Omega},\boldsymbol{K})$ is linearly sufficient, whereas both estimators are linearly complete.
Hence, the nonsingularity of $\boldsymbol{K}$ is closely related to the linear sufficiency of $\hat{\boldsymbol{\beta}}_{BL}(\boldsymbol{\Omega},\boldsymbol{K})$.\\
(ii) For a given $\boldsymbol{\Phi} \in \MS^+(n)$, $\hat{\boldsymbol{\beta}}_{BL}(\boldsymbol{\Phi},\boldsymbol{K})$
also possesses linear completeness.
However, $\hat{\boldsymbol{\beta}}_{BL}(\boldsymbol{\Phi},\boldsymbol{K})$ does not necessarily have linear sufficiency even if $\boldsymbol{K}\in \MS^+(k)$.
\end{rem}

\section{Identical equality between two Bayes linear estimators} \label{sec:4}
\subsection{Theoretical results}\label{sec:4.1}
In this subsection, we provide several necessary and sufficient conditions under which two BLEs coincide for all $\boldsymbol{y} \in \mathbb{R}^n$.

To explore the equality \eqref{B2E}, we begin with several technical lemmas.
First, we recall a useful representation of $\boldsymbol{\Omega}$ and $\boldsymbol{\Omega}^{-1}$.
For the proof, see \cite{RefK98} and \cite{RefKK04}.
\begin{lem}\label{Omega}
$\boldsymbol{\Omega} \in \MS^+(n)$ is written as
\begin{equation}\label{Orep}
 \boldsymbol{\Omega} 
= \left( \begin{matrix} \boldsymbol{X} & \boldsymbol{Z} \end{matrix}\right) \left(\begin{matrix} \boldsymbol{\Gamma} & \boldsymbol{\Xi} \\ \boldsymbol{\Xi}^\top & \boldsymbol{\Delta} \end{matrix} \right) \left( \begin{matrix} \boldsymbol{X}^\top \\ \boldsymbol{Z}^\top \end{matrix}\right)
\end{equation}
for some $\boldsymbol{\Gamma} \in \MS^+(k)$, $\boldsymbol{\Delta} \in \MS^+(n-k)$, and $k \times (n-k)$ matrix $\boldsymbol{\Xi}$.
Moreover, when $\boldsymbol{\Omega}$ is of the form \eqref{Orep}, it holds that
\[
\boldsymbol{\Omega}^{-1} 
= \left( \begin{matrix} \boldsymbol{X}(\boldsymbol{X}^\top \boldsymbol{X})^{-1} & \boldsymbol{Z}(\boldsymbol{Z}^\top \boldsymbol{Z})^{-1} \end{matrix}\right) \left(\begin{matrix} \boldsymbol{A} & \boldsymbol{B} \\ \boldsymbol{C} & \boldsymbol{D}\end{matrix} \right) \left( \begin{matrix} (\boldsymbol{X}^\top \boldsymbol{X})^{-1} \boldsymbol{X}^\top \\ (\boldsymbol{Z}^\top \boldsymbol{Z})^{-1} \boldsymbol{Z}^\top \end{matrix}\right),
\]
where
\begin{align*}
\boldsymbol{A} = (\boldsymbol{\Gamma} - \boldsymbol{\Xi} \boldsymbol{\Delta}^{-1} \boldsymbol{\Xi}^\top)^{-1}, \ 
\boldsymbol{B} = -\boldsymbol{A}\boldsymbol{\Xi} \boldsymbol{\Delta}^{-1}, \
\boldsymbol{C} = \boldsymbol{B}^\top, \
\boldsymbol{D} = \boldsymbol{\Delta}^{-1} + \boldsymbol{B}^\top \boldsymbol{A}^{-1} \boldsymbol{B}.
\end{align*}
\end{lem}

\begin{rem}
(i) The representation \eqref{Orep} holds when $\boldsymbol{\Gamma}, \boldsymbol{\Xi}, \boldsymbol{\Delta}$ are chosen as follows:
\begin{align*}
&\boldsymbol{\Gamma} = (\boldsymbol{X}^\top \boldsymbol{X})^{-1} \boldsymbol{X}^\top \boldsymbol{\Omega} \boldsymbol{X} (\boldsymbol{X}^\top \boldsymbol{X})^{-1}, \ 
\boldsymbol{\Xi} =  (\boldsymbol{X}^\top \boldsymbol{X})^{-1} \boldsymbol{X}^\top \boldsymbol{\Omega} \boldsymbol{Z} (\boldsymbol{Z}^\top \boldsymbol{Z})^{-1},\\
& \boldsymbol{\Delta} =  (\boldsymbol{Z}^\top \boldsymbol{Z})^{-1} \boldsymbol{Z}^\top \boldsymbol{\Omega} \boldsymbol{Z} (\boldsymbol{Z}^\top \boldsymbol{Z})^{-1}.
\end{align*}
(ii) The condition \eqref{Rcond} is equivalent to $\boldsymbol{\Xi} = \boldsymbol{0}$ in \eqref{Orep} and also to $\boldsymbol{X}^\top \boldsymbol{\Omega}\boldsymbol{Z} = \boldsymbol{0}$ and to $\boldsymbol{X}^\top \boldsymbol{\Omega}^{-1}\boldsymbol{Z} = \boldsymbol{0}$.\\
(iii) In this paper, we refer to $\boldsymbol{\Gamma},\boldsymbol{\Xi}, \boldsymbol{\Delta}$ as a covariance structure.
This covariance structure plays an important role in several results presented in Sections~\ref{sec:4} and \ref{sec:6}.
\end{rem}

Hereafter, we consider two models \eqref{model} introduced in Section~\ref{sec:1}.
From Lemma~\ref{Omega}, we obtain necessary and sufficient conditions for the equality \eqref{B2E}, which is one of the main results of this paper.

\begin{thm} \label{thm2}  
For given $\boldsymbol{K}_1, \boldsymbol{K}_2 \in \mathcal{S}^N(k)$, the equality \eqref{B2E} holds if and only if 
\begin{align}\label{B2Enew}
\boldsymbol{\Omega}\boldsymbol{X}\boldsymbol{K}_2 = \boldsymbol{X}\boldsymbol{K}_1,
\end{align}
which is also equivalent to
\begin{align}
&\boldsymbol{K}_2 \boldsymbol{X}^\top \boldsymbol{X} \boldsymbol{\Gamma} = \boldsymbol{K}_1 \label{414}\\
&\boldsymbol{K}_2 \boldsymbol{X}^\top \boldsymbol{X}\boldsymbol{\Xi}=\boldsymbol{0} \label{415}
\end{align}
in \eqref{Orep}.
\end{thm}
\begin{proof}
The equality \eqref{B2E} is equivalent to 
\begin{align}\label{4T}
& \boldsymbol{X}\hat{\boldsymbol{\beta}}_{BL}(\boldsymbol{\Omega}, \boldsymbol{K}_1) \equiv \boldsymbol{X}\hat{\boldsymbol{\beta}}_{BL}(\boldsymbol{I}_n, \boldsymbol{K}_2)\notag\\
\Leftrightarrow \quad & \boldsymbol{X}\boldsymbol{K}_1\boldsymbol{X}^{\top}(\boldsymbol{\Omega}+\boldsymbol{X}\boldsymbol{K}_1\boldsymbol{X}^{\top})^{-1}=\boldsymbol{X}\boldsymbol{K}_2\boldsymbol{X}^{\top}(\boldsymbol{I}_n+\boldsymbol{X}\boldsymbol{K}_2\boldsymbol{X}^{\top})^{-1}\notag\\
\Leftrightarrow \quad &\boldsymbol{\Omega}(\boldsymbol{\Omega}+\boldsymbol{X}\boldsymbol{K}_1\boldsymbol{X}^{\top})^{-1}=(\boldsymbol{I}_n+\boldsymbol{X}\boldsymbol{K}_2\boldsymbol{X}^{\top})^{-1}\notag\\
\Leftrightarrow \quad &(\boldsymbol{I}_n+\boldsymbol{X}\boldsymbol{K}_2\boldsymbol{X}^{\top})\boldsymbol{\Omega}=\boldsymbol{\Omega}+\boldsymbol{X}\boldsymbol{K}_1\boldsymbol{X}^{\top}\notag\\
\Leftrightarrow \quad &\boldsymbol{X}\boldsymbol{K}_2\boldsymbol{X}^{\top}\boldsymbol{\Omega}=\boldsymbol{X}\boldsymbol{K}_1\boldsymbol{X}^{\top}\notag\\
\Leftrightarrow \quad &\boldsymbol{K}_2\boldsymbol{X}^{\top}\boldsymbol{\Omega}=\boldsymbol{K}_1\boldsymbol{X}^{\top}.
\end{align}
Taking transposes yields \eqref{B2Enew}.
Since $\boldsymbol{X}^\top \boldsymbol{Z} = \boldsymbol{0}$ and the matrix $(\begin{matrix}
\boldsymbol{X} & \boldsymbol{Z}
\end{matrix})$ is nonsingular, the equality \eqref{4T} holds if and only if both of the following equalities hold:
\begin{align}
& \boldsymbol{X}^\top  \boldsymbol{X} \boldsymbol{K}_2 \boldsymbol{X}^\top \boldsymbol{\Omega} \boldsymbol{X}=\boldsymbol{X}^\top \boldsymbol{X} \boldsymbol{K}_1 \boldsymbol{X}^\top \boldsymbol{X}, \label{XX}\\
&  \boldsymbol{X}^\top  \boldsymbol{X} \boldsymbol{K}_2 \boldsymbol{X}^\top \boldsymbol{\Omega} \boldsymbol{Z} = \boldsymbol{0}. \label{XZ}
\end{align}
From Lemma~\ref{Omega}, the conditions \eqref{XX} and \eqref{XZ} are rewritten as 
\begin{align*}
&\boldsymbol{X}^\top \boldsymbol{X} \boldsymbol{K}_2 \boldsymbol{X}^\top \boldsymbol{X} \boldsymbol{\Gamma} \boldsymbol{X}^\top \boldsymbol{X} = \boldsymbol{X}^\top \boldsymbol{X} \boldsymbol{K}_1 \boldsymbol{X}^\top  \boldsymbol{X}\\
\Leftrightarrow \quad & \boldsymbol{K}_2 \boldsymbol{X}^\top \boldsymbol{X} \boldsymbol{\Gamma} = \boldsymbol{K}_1
\end{align*}
and 
\begin{align*}
&\boldsymbol{X}^\top \boldsymbol{X} \boldsymbol{K}_2 \boldsymbol{X}^\top \boldsymbol{X} \boldsymbol{\Xi} \boldsymbol{Z}^\top \boldsymbol{Z} = \boldsymbol{0}\\
\Leftrightarrow \quad & \boldsymbol{K}_2 \boldsymbol{X}^\top \boldsymbol{X} \boldsymbol{\Xi}  = \boldsymbol{0}.
\end{align*}
This completes the proof.
\end{proof}

\begin{rem}
(i) Theorem~\ref{thm2} shows that it is sufficient to identify the covariance structure $\boldsymbol{\Gamma}$ and $\boldsymbol{\Xi}$ to check the equality \eqref{B2E}, even if $\boldsymbol{\Delta}$ remains unknown. 
See also Example~\ref{ex1} in Section~\ref{sec:4.2}.\\
(ii) The equality \eqref{B2E} is also equivalent to $\hat{\boldsymbol{\beta}}_{BL}(\boldsymbol{\Omega}, \boldsymbol{K}_1) = \Hat{\boldsymbol{\beta}}_{BL}(\boldsymbol{I}_n, \boldsymbol{K}_2)$ almost surely.
\end{rem}

If $\boldsymbol{K}_1,\boldsymbol{K}_2 \in \MS^+(k)$, we directly obtain the following corollary to Theorem~\ref{thm2}.
\begin{cor}\label{cor34}
(i) For given $\boldsymbol{K}_1,\boldsymbol{K}_2 \in \MS^+(k)$, the equality \eqref{B2E} holds if and only if $\boldsymbol{X}=\boldsymbol{\Omega X }\boldsymbol{K}_2\boldsymbol{K}_1^{-1}$,
which is also equivalent to 
\begin{align*}
\boldsymbol{\Omega}=\boldsymbol{X}(\boldsymbol{X}^\top \boldsymbol{X})^{-1}\boldsymbol{K}_2^{-1}\boldsymbol{K}_1 \boldsymbol{X}^\top + \boldsymbol{Z}\boldsymbol{\Delta}\boldsymbol{Z}^\top \quad \fs \ \ \boldsymbol{\Delta} \in \MS^+(n-k).    
\end{align*}
(ii) For a given $\boldsymbol{K} \in \MS^+(k)$, the equality 
$\hat{\boldsymbol{\beta}}_{BL}(\boldsymbol{\Omega}, \boldsymbol{K}) \equiv \hat{\boldsymbol{\beta}}_{BL} (\boldsymbol{I}_n, \boldsymbol{K})$
holds if and only if $\boldsymbol{X} = \boldsymbol{\Omega X}$, which is also equivalent to 
\begin{align*}
\boldsymbol{\Omega}= \boldsymbol{P}_X  + \boldsymbol{Z}\boldsymbol{\Delta}\boldsymbol{Z}^\top \quad \fs \ \ \boldsymbol{\Delta} \in \MS^+(n-k),  
\end{align*}
where 
\begin{align}\label{P}
\boldsymbol{P}_X = \boldsymbol{X}(\boldsymbol{X}^\top \boldsymbol{X})^{-1}\boldsymbol{X}^\top.
\end{align}
\end{cor}

\begin{rem}
$\boldsymbol{X} = \boldsymbol{\Omega X}$ is further equivalent to $\boldsymbol{P}_X = \boldsymbol{P}_X \boldsymbol{\Omega}$.
\end{rem}

\begin{rem}[\textbf{Equality between ordinary ridge estimators}]
Let $\boldsymbol{K}=\lambda \boldsymbol{I}_k \in \MS^+(k)$. 
Then Corollary~\ref{cor34} implies that $\hat{\boldsymbol{\beta}}_{BL}(\boldsymbol{\Omega},\lambda \boldsymbol{I}_k) \equiv \hat{\boldsymbol{\beta}}_{BL}(\boldsymbol{I}_n, \lambda \boldsymbol{I}_k)$ 
holds if and only if $\boldsymbol{X} = \boldsymbol{\Omega X}$, which is also
equivalent to $\boldsymbol{\Omega}= \boldsymbol{P}_X  + \boldsymbol{Z}\boldsymbol{\Delta}\boldsymbol{Z}^\top$ for some $\boldsymbol{\Delta} \in \MS^+(n-k)$.
\end{rem}

\begin{rem}[\textbf{Equality between typical shrinkage estimators}]
\mbox{}\\
Let $\boldsymbol{K}_1 = \rho(\boldsymbol{X}^\top \boldsymbol{\Omega}^{-1}\boldsymbol{X})^{-1}$ and $\boldsymbol{K}_2 = \rho(\boldsymbol{X}^\top \boldsymbol{X})^{-1}$, where $\rho$ is a positive constant.
Suppose that $\boldsymbol{\Omega}$ is of the form \eqref{Rcond}.
Then the two matrices satisfy \eqref{414} and \eqref{415}.
In fact, it follows from Lemma~\ref{Omega} that $\boldsymbol{X}^\top \boldsymbol{\Omega}^{-1}\boldsymbol{X} = \boldsymbol{\Gamma}^{-1}$ and thus $\boldsymbol{K}_1 = \rho \boldsymbol{\Gamma}$.
So, it holds that $\boldsymbol{K}_2\boldsymbol{X}^\top \boldsymbol{X}\boldsymbol{\Gamma} = \rho(\boldsymbol{X}^\top \boldsymbol{X})^{-1}\boldsymbol{X}^\top \boldsymbol{X}\boldsymbol{\Gamma} = \rho\boldsymbol{\Gamma} = \boldsymbol{K}_1$.
Moreover, \eqref{Rcond} is equivalent to $\boldsymbol{\Xi} = \boldsymbol{0}$, implying that \eqref{415} holds.
Hence, the equality $\hat{\boldsymbol{\beta}}_{BL}(\boldsymbol{\Omega}, \rho(\boldsymbol{X}^\top \boldsymbol{\Omega}^{-1}\boldsymbol{X})^{-1}) = \hat{\boldsymbol{\beta}}_{BL}(\boldsymbol{I}_n, \rho(\boldsymbol{X}^\top \boldsymbol{X})^{-1})$ is satisfied.
More specifically, \eqref{Rcond} is the necessary and sufficient condition for $\hat{\boldsymbol{\beta}}_{BL}(\boldsymbol{\Omega}, \rho(\boldsymbol{X}^\top \boldsymbol{\Omega}^{-1}\boldsymbol{X})^{-1}) \equiv \hat{\boldsymbol{\beta}}_{BL}(\boldsymbol{I}_n, \rho(\boldsymbol{X}^\top \boldsymbol{X})^{-1})$.
This conclusion follows directly from the forms of the shrinkage estimators.
\end{rem}

Hereafter, we show that the equality \eqref{B2E} is closely related to Bayes risk efficiency and to the expectation vector and covariance matrices of BLEs.
By analogy with the usual notion of efficiency, we define the Bayes risk efficiency by using the trace of the Bayes risk matrix:
\begin{defn}
The Bayes risk efficiency of $\delta_1$ relative to $\delta_2$ with respect to $\pi$ is defined as
\begin{align*}
\BRE(\delta_1;\delta_2) = \frac{\tr\left(\BRM(\delta_2)\right)}{\tr\left(\BRM(\delta_1)\right)} = \cfrac{\E_\pi[\E[\|\delta_2 - \boldsymbol{\beta}\|^2 \ | \ \boldsymbol{\beta}, \sigma^2]]}{\E_\pi[\E[\|\delta_1 - \boldsymbol{\beta}\|^2 \ | \ \boldsymbol{\beta}, \sigma^2]]},
\end{align*}
where $\tr(\boldsymbol{M})$ denotes the trace of a matrix $\boldsymbol{M}$.
\end{defn}
It is obvious that for given $\boldsymbol{K}_1, \boldsymbol{K}_2 \in \MS^N(k)\setminus{\{\boldsymbol{0}\}}$, 
\begin{align*}
0 < \BRE(\hat{\boldsymbol{\beta}}_{BL}(\boldsymbol{I}_n, \boldsymbol{K}_2);\hat{\boldsymbol{\beta}}_{BL}(\boldsymbol{\Omega}, \boldsymbol{K}_1)) \leq 1
\end{align*}
holds under Model 1.
In particular, if $\boldsymbol{W}_1,\boldsymbol{W}_2\in\MS^+(k)$ are fixed and
$\gamma_1, \ \gamma_2 \downarrow 0$, then the Bayes risk efficiency reduces to the risk efficiency of OLSE relative to GLSE based on the  mean squared error.
Proposition~\ref{pr45} shows that the equality \eqref{B2E} is closely related to the maximization of Bayes risk efficiency.
\begin{prop}\label{pr45}
For given $\boldsymbol{K}_1,\boldsymbol{K}_2 \in \MS^N(k)\setminus{\{\boldsymbol{0}\}}$, the equality \eqref{B2E} is equivalent to 
\begin{align}\label{MBR1}
\BRE(\hat{\boldsymbol{\beta}}_{BL}(\boldsymbol{I}_n, \boldsymbol{K}_2);\hat{\boldsymbol{\beta}}_{BL}(\boldsymbol{\Omega}, \boldsymbol{K}_1)) = 1
\end{align}
under Model 1.
\end{prop}
\begin{proof}
To simplify the notation, let us denote
\begin{align*}
\boldsymbol{L}_\Omega = \boldsymbol{K}_1\boldsymbol{X}^\top(\boldsymbol{\Omega} + \boldsymbol{X}\boldsymbol{K}_1\boldsymbol{X}^\top)^{-1}, \quad \boldsymbol{L}_I = \boldsymbol{K}_2\boldsymbol{X}^\top(\boldsymbol{I}_n + \boldsymbol{X}\boldsymbol{K}_2\boldsymbol{X}^\top)^{-1},
\end{align*}
and $\boldsymbol{T}_1 = \boldsymbol{\Omega} + \boldsymbol{X}\boldsymbol{K}_1\boldsymbol{X}^\top$.
It follows from \eqref{BRM2} in Proposition~\ref{pr21} that \eqref{MBR1} is equivalent to
\begin{align*}
&\tr\left[\BRM(\hat{\boldsymbol{\beta}}_{BL}(\boldsymbol{\Omega}, \boldsymbol{K}_1))\right] = \tr\left[\BRM(\hat{\boldsymbol{\beta}}_{BL}(\boldsymbol{I}_n, \boldsymbol{K}_2)\right]\\
\Leftrightarrow \quad & \tr\left[\BRM(\hat{\boldsymbol{\beta}}_{BL}(\boldsymbol{\Omega}, \boldsymbol{K}_1))\right] =  \gamma_1\tr\left[(\boldsymbol{L}_I - \boldsymbol{L}_\Omega) \boldsymbol{T}_1 (\boldsymbol{L}_I - \boldsymbol{L}_\Omega)^\top \right] + \tr\left[\BRM(\hat{\boldsymbol{\beta}}_{BL}(\boldsymbol{\Omega}, \boldsymbol{K}_1))\right]\\
\Leftrightarrow \quad & \tr\left[(\boldsymbol{L}_I - \boldsymbol{L}_\Omega) \boldsymbol{T}_1 (\boldsymbol{L}_I - \boldsymbol{L}_\Omega)^\top \right] = 0\\
\Leftrightarrow \quad & \boldsymbol{L}_I - \boldsymbol{L}_\Omega = \boldsymbol{0}\\
\Leftrightarrow \quad & \hat{\boldsymbol{\beta}}_{BL}(\boldsymbol{\Omega}, \boldsymbol{K}_1) \equiv \hat{\boldsymbol{\beta}}_{BL}(\boldsymbol{I}_n, \boldsymbol{K}_2),
\end{align*}
where $\BRM(\boldsymbol{L}\boldsymbol{y}) = \gamma_1\{\boldsymbol{L}\boldsymbol{\Omega}\boldsymbol{L}^\top + (\boldsymbol{L}\boldsymbol{X} - \boldsymbol{I}_k) \boldsymbol{K}_1(\boldsymbol{L}\boldsymbol{X} - \boldsymbol{I}_k)^\top\}$.
This completes the proof.
\end{proof}

\begin{rem}
Proposition~\ref{pr45} reveals that there is no loss of Bayes risk efficiency in using $\hat{\boldsymbol{\beta}}_{BL}(\boldsymbol{I}_n,\boldsymbol{K}_2)$ under Model 1 if and only if it coincides with $\hat{\boldsymbol{\beta}}_{BL}(\boldsymbol{\Omega},\boldsymbol{K}_1)$ for all $\boldsymbol{y}\in\mathbb{R}^n$.
\end{rem}

Next, we provide a moment-based characterization of the equality \eqref{B2E}.
\begin{prop}\label{pr47}
For given $\boldsymbol{K}_1,\boldsymbol{K}_2 \in \MS^N(k)$, the equality \eqref{B2E} holds if and only if
\begin{align}
& \E\left[\hat{\boldsymbol{\beta}}_{BL}(\boldsymbol{\Omega}, \boldsymbol{K}_1)\right] = \E\left[\hat{\boldsymbol{\beta}}_{BL}(\boldsymbol{I}_n, \boldsymbol{K}_2)\right] \quad \fa  \ \ \boldsymbol{\beta} \in \mathbb{R}^k, \label{Eeq}\\
& \V\left(\hat{\boldsymbol{\beta}}_{BL}(\boldsymbol{\Omega}, \boldsymbol{K}_1)\right) = \V\left(\hat{\boldsymbol{\beta}}_{BL}(\boldsymbol{I}_n, \boldsymbol{K}_2)\right), \label{Ceq}
\end{align}
where both moments are calculated under Model~1.
\end{prop}
\begin{proof}
To simplify the notation, let us denote
\begin{align*}
\boldsymbol{S}_1 = \boldsymbol{K}_1\left[(\boldsymbol{X}^\top \boldsymbol{\Omega}^{-1}\boldsymbol{X})^{-1} + \boldsymbol{K}_1\right]^{-1}, \quad \boldsymbol{S}_2 = \boldsymbol{K}_2\left[(\boldsymbol{X}^\top \boldsymbol{X})^{-1} + \boldsymbol{K}_2\right]^{-1}
\end{align*}
From Lemma~\ref{lem24} and Remark~\ref{rm3}-(ii), \eqref{Eeq} can be rewritten as
\begin{align}\label{Eeq2}
\boldsymbol{S}_1 = \boldsymbol{S}_2.
\end{align}
Under \eqref{Eeq2}, \eqref{Ceq} is equivalent to
\begin{align}\label{Ceq2}
&\boldsymbol{S}_1 (\boldsymbol{X}^\top \boldsymbol{\Omega}^{-1}\boldsymbol{X})^{-1} \boldsymbol{S}_1^\top = \boldsymbol{S}_2 (\boldsymbol{X}^\top \boldsymbol{X})^{-1} \boldsymbol{X}^\top \boldsymbol{\Omega}\boldsymbol{X}(\boldsymbol{X}^\top \boldsymbol{X})^{-1}\boldsymbol{S}_2^\top \notag\\
\Leftrightarrow \quad & \boldsymbol{S}_1 (\boldsymbol{\Gamma} - \boldsymbol{\Xi}\boldsymbol{\Delta}^{-1}\boldsymbol{\Xi}^\top)\boldsymbol{S}_1^\top = \boldsymbol{S}_1 \boldsymbol{\Gamma}\boldsymbol{S}_1^\top \notag\\
\Leftrightarrow \quad & \boldsymbol{S}_1 \boldsymbol{\Xi}\boldsymbol{\Delta}^{-1}\boldsymbol{\Xi}^\top \boldsymbol{S}_1^\top = \boldsymbol{0},
\end{align}
where $\boldsymbol{\Gamma}, \boldsymbol{\Xi}$, and $\boldsymbol{\Delta}$ are defined in \eqref{Orep}.
Equivalently, \eqref{Ceq2} can be written as 
\begin{align}\label{Ceq3}
\boldsymbol{S}_1\boldsymbol{\Xi} = \boldsymbol{0}.
\end{align}
On the other hand, the equality \eqref{B2E} is equivalent to
\begin{align*}
&\boldsymbol{S}_1 (\boldsymbol{X}^\top \boldsymbol{\Omega}^{-1}\boldsymbol{X})^{-1} \boldsymbol{X}^\top \boldsymbol{\Omega}^{-1} = \boldsymbol{S}_2(\boldsymbol{X}^\top \boldsymbol{X})^{-1} \boldsymbol{X}^\top \\
\Leftrightarrow \quad & \boldsymbol{S}_1 (\boldsymbol{X}^\top \boldsymbol{\Omega}^{-1}\boldsymbol{X})^{-1} \boldsymbol{X}^\top \boldsymbol{\Omega}^{-1}
\begin{pmatrix}
\boldsymbol{X} & \boldsymbol{Z}
\end{pmatrix} 
= \boldsymbol{S}_2(\boldsymbol{X}^\top \boldsymbol{X})^{-1} \boldsymbol{X}^\top 
\begin{pmatrix}
\boldsymbol{X} & \boldsymbol{Z}
\end{pmatrix}.
\end{align*}
Noting that $(\boldsymbol{X}^\top \boldsymbol{\Omega}^{-1}\boldsymbol{X})^{-1}\boldsymbol{X}^\top \boldsymbol{\Omega}^{-1}\boldsymbol{Z} = -\boldsymbol{\Xi}\boldsymbol{\Delta}^{-1}$ and $\boldsymbol{X}^\top \boldsymbol{Z} = \boldsymbol{0}$,
the preceding block-matrix equality yields \eqref{Eeq2} and \eqref{Ceq3}, respectively.
This completes the proof.
\end{proof}

\begin{rem}
Since OLSE and GLSE are unbiased for $\boldsymbol{\beta}$, the equality of these estimators is characterized by the equality of covariance matrices; see, for example, \cite{RefPS89} and \cite{RefP05}. 
Proposition~\ref{pr47} extends this result to biased BLEs by additionally requiring the equality of their expectation vectors.
\end{rem}

The following corollary provides a necessary and sufficient condition for the equality between two classes of BLEs from Theorem~\ref{thm2}.
\begin{cor}\label{cor35}
\mbox{}\\
Consider the classes $\{\hat{\boldsymbol{\beta}}_{BL}(\boldsymbol{\Omega}, \boldsymbol{K})\}_{\boldsymbol{K}\in \mathcal{S}^{N}(k)}$ and $\{\hat{\boldsymbol{\beta}}_{BL}(\boldsymbol{I}_n, \boldsymbol{K})\}_{\boldsymbol{K}\in \mathcal{S}^{N}(k)}$ of BLEs.
The equality 
\begin{align}\label{CE}
\{\hat{\boldsymbol{\beta}}_{BL}(\boldsymbol{\Omega}, \boldsymbol{K})\}_{\boldsymbol{K}\in \mathcal{S}^{N}(k)} = \{\hat{\boldsymbol{\beta}}_{BL}(\boldsymbol{I}_n, \boldsymbol{K})\}_{\boldsymbol{K}\in \mathcal{S}^{N}(k)}
\end{align}
holds if and only if $\boldsymbol{\Omega}$ is of the form 
\begin{align}\label{Orep2}
\boldsymbol{\Omega} = c \boldsymbol{P}_X  + \boldsymbol{Z}\boldsymbol{\Delta}\boldsymbol{Z}^\top
\end{align}
for some $c > 0$ and $\boldsymbol{\Delta}\in \MS^+(n-k)$, where $\boldsymbol{P}_X$ is defined in \eqref{P}.
\end{cor}
\begin{proof}
To see the necessity, suppose that the equality \eqref{CE} is satisfied.
Let $\boldsymbol{M} = \boldsymbol{X}^\top \boldsymbol{X}\boldsymbol{\Gamma}$.
Since $\hat{\boldsymbol{\beta}}_{BL}(\boldsymbol{I}_n, \boldsymbol{K}) \in \{\hat{\boldsymbol{\beta}}_{BL}(\boldsymbol{\Omega}, \boldsymbol{K})\}_{\boldsymbol{K}\in \mathcal{S}^{N}(k)}$, there exists a matrix $\Tilde{\boldsymbol{K}} \in \mathcal{S}^{N}(k)$ such that
$\hat{\boldsymbol{\beta}}_{BL}(\boldsymbol{\Omega}, \Tilde{\boldsymbol{K}}) = \hat{\boldsymbol{\beta}}_{BL}(\boldsymbol{I}_n, \boldsymbol{K})$ for each $\boldsymbol{K}\in\MS^N(k)$, which is further written from Theorem~\ref{thm2} as follows:
\begin{align}
&\boldsymbol{K}\boldsymbol{M} = \Tilde{\boldsymbol{K}}, \label{cond1}\\
&\boldsymbol{K}\boldsymbol{X}^\top \boldsymbol{X}\boldsymbol{\Xi} = \boldsymbol{0}. \label{cond2}
\end{align}
If $\boldsymbol{K} = \boldsymbol{I}_k$, it follows from \eqref{cond1} and \eqref{cond2} that $\boldsymbol{M} \in \MS^+(k)$ and $\boldsymbol{\Xi} = \boldsymbol{0}$.
Hereafter, we show that $\boldsymbol{\Gamma} = c(\boldsymbol{X}^\top \boldsymbol{X})^{-1}$.
Let $\boldsymbol{K} = \boldsymbol{e}_i\boldsymbol{e}_i^\top$ for $i = 1, \ldots, k$, where $\boldsymbol{e}_i$ denotes the $i$th standard basis vector in $\mathbb{R}^k$.
Since $\boldsymbol{K}\boldsymbol{M} \in \MS^N(k)$, we have
\begin{align*}
\boldsymbol{K}\boldsymbol{M} = (\boldsymbol{K}\boldsymbol{M})^\top = \boldsymbol{M}\boldsymbol{K} \quad \Leftrightarrow \quad \boldsymbol{e}_i \boldsymbol{e}_i^\top \boldsymbol{M} = \boldsymbol{M} \boldsymbol{e}_i \boldsymbol{e}_i^\top.
\end{align*}
Postmultiplying by $\boldsymbol{e}_i$ yields
\begin{align*}
\boldsymbol{e}_i \boldsymbol{e}_i^\top \boldsymbol{M}\boldsymbol{e}_i  = \boldsymbol{M} \boldsymbol{e}_i \boldsymbol{e}_i^\top \boldsymbol{e}_i 
 \quad \Leftrightarrow \quad \boldsymbol{M}\boldsymbol{e}_i = \boldsymbol{e}_i (\boldsymbol{e}_i^\top \boldsymbol{M}\boldsymbol{e}_i )
\quad \Leftrightarrow \quad  \boldsymbol{M}\boldsymbol{e}_i =  m_i \boldsymbol{e}_i \quad (\say).
\end{align*}
Therefore, we obtain
\begin{align*}
\boldsymbol{M} = (\boldsymbol{M}\boldsymbol{e}_1, \ldots, \boldsymbol{M}\boldsymbol{e}_k)  = (m_1 \boldsymbol{e}_1, \ldots, m_k \boldsymbol{e}_k) = \diag(m_1, \ldots, m_k),
\end{align*}
where $\diag(\cdot)$ denotes a diagonal matrix.
Also, we consider $\boldsymbol{K} = (\boldsymbol{e}_i + \boldsymbol{e}_j)(\boldsymbol{e}_i + \boldsymbol{e}_j)^\top$ for any different $ 1\leq i,j \leq k$.
Note that $\boldsymbol{K}$ can be written as
\begin{align*}
    \boldsymbol{K} = \boldsymbol{e}_i\boldsymbol{e}_i^\top + \boldsymbol{e}_i \boldsymbol{e}_j^\top + \boldsymbol{e}_j \boldsymbol{e}_i^\top + \boldsymbol{e}_j\boldsymbol{e}_j^\top,
\end{align*}
it is sufficient to compare the $(i,j)$th entries of $\boldsymbol{K}\boldsymbol{M}$ and $\boldsymbol{M}\boldsymbol{K}$.
These entries are $m_i$ and $m_j$, respectively, and therefore $m_i=m_j$.
This implies that there exists a positive constant such that $m_i = c$ for all $i = 1, \ldots, k$.
Hence, we have $\boldsymbol{M} = c\boldsymbol{I}_k \ \Leftrightarrow \ \boldsymbol{X}^\top \boldsymbol{X} \boldsymbol{\Gamma} = c\boldsymbol{I}_k \ \Leftrightarrow \ \boldsymbol{\Gamma} = c(\boldsymbol{X}^\top \boldsymbol{X})^{-1}$.

To see the sufficiency, suppose that $\boldsymbol{\Omega}$ is of the form \eqref{Orep2}.
Theorem~\ref{thm2} implies that the equality
$\hat{\boldsymbol{\beta}}_{BL}(\boldsymbol{I}_n, \boldsymbol{K}) = \hat{\boldsymbol{\beta}}_{BL}(\boldsymbol{\Omega}, c\boldsymbol{K})$ holds for any $\boldsymbol{K}$, and hence
\begin{align*}
\{\hat{\boldsymbol{\beta}}_{BL}(\boldsymbol{I}_n, \boldsymbol{K})\}_{\boldsymbol{K}\in \mathcal{S}^{N}(k)} \subseteq \{\hat{\boldsymbol{\beta}}_{BL}(\boldsymbol{\Omega}, \boldsymbol{K})\}_{\boldsymbol{K}\in \MS^N(k)}.
\end{align*}
Since the  equality
$\hat{\boldsymbol{\beta}}_{BL}(\boldsymbol{\Omega},\boldsymbol{K}) = \hat{\boldsymbol{\beta}}_{BL}(\boldsymbol{I}_n, c^{-1}\boldsymbol{K})$
holds, the reverse inclusion is true.
Therefore, the equality \eqref{CE} holds.
This completes the proof.
\end{proof}

\begin{rem}
\mbox{}\\
When $\boldsymbol{\Omega}$ is of the form \eqref{Orep2}, it suffices to consider the class
$\{\hat{\boldsymbol{\beta}}_{BL}(\boldsymbol{I}_n, \boldsymbol{K})\}_{\boldsymbol{K}\in \mathcal{S}^{N}(k)}$.
\end{rem}
\subsection{Applications}\label{sec:4.2}
In this subsection, as applications of Theorem~\ref{thm2}, we present concrete examples in which the equality \eqref{B2E} holds.
Hereafter, we consider the case where $\boldsymbol{\Omega}$ is partially unknown.
In this case, we often adopt a two-step estimation for $\boldsymbol{\beta}$.
However, when the equality \eqref{B2E} holds, we can dispense with this procedure.

First, we present a numerical example in the simple case $n =3$ and  $k=2$.
\begin{exam}[\textbf{Numerical example for Theorem~\ref{thm2}}]\label{ex1}
Let  
\begin{align*}
\boldsymbol{X}=
\begin{pmatrix}
1 & 0 \\
0 & 1 \\
0 & 0 
\end{pmatrix},
\boldsymbol{Z} =
\begin{pmatrix}
0\\
0\\
2
\end{pmatrix}, \
\boldsymbol{\Omega} =
\begin{pmatrix}
2 & 1 & 4 \\
1 & 2 & -4 \\
4 & -4 & 4a 
\end{pmatrix}, \ 
\boldsymbol{K}_1=
\begin{pmatrix}
3 & 3 \\
3 & 3 
\end{pmatrix}, \
\boldsymbol{K}_2= \begin{pmatrix}
1 & 1 \\
1 & 1 
\end{pmatrix},
\end{align*}
where $a > 8$ is an unknown constant.
When we choose the covariance structure 
\begin{align*}
\boldsymbol{\Gamma} =
\begin{pmatrix}
2 & 1 \\
1 & 2
\end{pmatrix}, \
\boldsymbol{\Xi} =
\begin{pmatrix}
2\\
-2
\end{pmatrix}, \ \boldsymbol{\Delta} = a,
\end{align*}
$\boldsymbol{\Omega}$ can be written in the form of \eqref{Orep}. Since $\boldsymbol{X}^\top\boldsymbol{X}=\boldsymbol{I}_2$, we have
\begin{align*}
\boldsymbol{K}_2\boldsymbol{X}^\top\boldsymbol{X}\boldsymbol{\Gamma} 
= \begin{pmatrix}
1 & 1 \\
1 & 1 
\end{pmatrix}\begin{pmatrix}
2 & 1 \\
1 & 2
\end{pmatrix}
= \begin{pmatrix}
3 & 3 \\
3 & 3 
\end{pmatrix}
=\boldsymbol{K}_1
\end{align*}
and 
\begin{align*}
\boldsymbol{K}_2 \boldsymbol{X}^\top\boldsymbol{X}\boldsymbol{\Xi} 
=
\begin{pmatrix}
1 & 1 \\
1 & 1 
\end{pmatrix}
\begin{pmatrix}
2\\
-2
\end{pmatrix}
= \begin{pmatrix}
0\\
0
\end{pmatrix}.
\end{align*}
This implies that both \eqref{414} and \eqref{415} are satisfied.
Hence, from Theorem~\ref{thm2}, the following  equality holds:
\begin{align*}
\hat{\boldsymbol{\beta}}_{BL}(\boldsymbol{\Omega}, \boldsymbol{K}_1) \equiv \hat{\boldsymbol{\beta}}_{BL} (\boldsymbol{I}_3, \boldsymbol{K}_2)\left(=\cfrac{1}{3}
\begin{pmatrix}
1 & 1 & 0\\
1 & 1 & 0
\end{pmatrix}\boldsymbol{y}\right).
\end{align*}
\end{exam}
Although $\boldsymbol{\Omega}$ contains the unknown parameter $a$ through $\boldsymbol{\Delta}=a$, this parameter is irrelevant to the verification of \eqref{414} and \eqref{415}. 
So, we can avoid the estimation of $\boldsymbol{\Omega}$.

Next, we impose specific structures on $\boldsymbol{\Omega}$ in  \eqref{glm} and investigate conditions under which the equality 
\begin{align}\label{B2EK}
\hat{\boldsymbol{\beta}}_{BL}(\hat{\boldsymbol{\Omega}},\boldsymbol{K}) \equiv \hat{\boldsymbol{\beta}}_{BL}(\boldsymbol{I}_n,\boldsymbol{K})
\end{align}
holds.
Note that $\hat{\boldsymbol{\Omega}} = \boldsymbol{\Omega}(\hat{\omega})$, where $\hat{\omega}$ is an estimator of an unknown parameter $\omega$.
\begin{exam}[\textbf{Rao's mixed-effects model \citep{RefR67}}]\label{ex2}
\mbox{}\\
Consider Rao's mixed-effects model
\begin{align}\label{Raom}
\boldsymbol{y} = \boldsymbol{X}\boldsymbol{\beta} + \boldsymbol{X}\boldsymbol{\gamma} + \boldsymbol{Z}\boldsymbol{\delta} + \boldsymbol{\ve},
\end{align}
where $\boldsymbol{\gamma},\boldsymbol{\delta}$, and $\boldsymbol{\ve}$ are mutually uncorrelated random vectors with $\E[\boldsymbol{\gamma}] = \boldsymbol{0}_k,\ \V(\boldsymbol{\gamma}) = \sigma^2\Bar{\boldsymbol{\Gamma}} \in \MS^N(k), \ \E[\boldsymbol{\delta}] = \boldsymbol{0}_{n-k}, \ \V(\boldsymbol{\delta}) = \sigma^2\Bar{\boldsymbol{\Delta}}\in \MS^N(n-k), \ \E[\boldsymbol{\ve}] = \boldsymbol{0}_n, \ \V(\boldsymbol{\ve}) = \sigma^2 \boldsymbol{I}_n$. 
Suppose that $\Bar{\boldsymbol{\Gamma}}$ is known, whereas $\Bar{\boldsymbol{\Delta}}$ is unknown.
Then $\boldsymbol{\Omega}$ can be written as
\begin{align*}
\boldsymbol{\Omega} = \boldsymbol{I}_n + \boldsymbol{X}\Bar{\boldsymbol{\Gamma}}\boldsymbol{X}^\top + \boldsymbol{Z}\Bar{\boldsymbol{\Delta}}\boldsymbol{Z}^\top 
= \boldsymbol{X}\left[(\boldsymbol{X}^\top \boldsymbol{X})^{-1} + \Bar{\boldsymbol{\Gamma}}\right]\boldsymbol{X}^\top + \boldsymbol{Z}\left[(\boldsymbol{Z}^\top \boldsymbol{Z})^{-1} + \Bar{\boldsymbol{\Delta}}\right]\boldsymbol{Z}^\top,
\end{align*}
which implies that $\boldsymbol{\Omega}$ is of the form \eqref{Rcond} and hence \eqref{415} is satisfied.
Moreover, by imposing an additional assumption on $\Bar{\boldsymbol{\Gamma}}$, we obtain the following corollary.
\begin{cor}\label{cor49}
The equality \eqref{B2EK} holds under \eqref{Raom} if and only if $\boldsymbol{K}\boldsymbol{X}^\top\boldsymbol{X}\Bar{\boldsymbol{\Gamma}} = \boldsymbol{0}$.
\end{cor}
\begin{proof}
Since
\begin{align*}
\boldsymbol{\Omega}\boldsymbol{X}\boldsymbol{K} &= \boldsymbol{X}\left[(\boldsymbol{X}^\top \boldsymbol{X})^{-1} + \Bar{\boldsymbol{\Gamma}}\right]\boldsymbol{X}^\top \boldsymbol{X}\boldsymbol{K}\\
& = \boldsymbol{X}\boldsymbol{K} +\boldsymbol{X}\Bar{\boldsymbol{\Gamma}}\boldsymbol{X}^\top\boldsymbol{X}\boldsymbol{K},
\end{align*}
$\boldsymbol{\Omega}\boldsymbol{X}\boldsymbol{K} = \boldsymbol{X}\boldsymbol{K}$ is equivalent to $\boldsymbol{K}\boldsymbol{X}^\top\boldsymbol{X}\Bar{\boldsymbol{\Gamma}} = \boldsymbol{0}$.
Hence, it follows from Theorem~\ref{thm2} that  the equality \eqref{B2EK} holds.
This completes the proof.
\end{proof}
\begin{rem}
(i) From Corollary~\ref{cor49}, we can check whether the equality \eqref{B2EK} holds even when $\Bar{\boldsymbol{\Delta}}$ is unknown.\\
(ii) Suppose that $\boldsymbol{K} \in \MS^+(k)$.
Then the equality \eqref{B2EK} holds if and only if $\Bar{\boldsymbol{\Gamma}} = \boldsymbol{0}$.
Furthermore, when we consider the model 
\begin{align}\label{LMM}
\boldsymbol{y} = \boldsymbol{X}\boldsymbol{\beta} + \boldsymbol{Z}\boldsymbol{\delta} + \boldsymbol{\ve},
\end{align}
the equality \eqref{B2EK} is always satisfied for $\boldsymbol{K} \in \MS^N(k)$.
\end{rem}
\textbf{Numerical Experiment.}
We numerically examine the equivalence \eqref{B2EK} under \eqref{LMM}.
Let $n=20$ and $k=3$, and consider the quadratic-trend design matrix
\begin{align*}
\boldsymbol{X}
=
\begin{pmatrix}
1 & t_1 & t_1^2\\
\vdots & \vdots & \vdots\\
1 & t_{20} & t_{20}^2
\end{pmatrix},
\qquad
t_i=-1+\frac{2(i-1)}{19},
\quad i=1,\ldots,20.
\end{align*}
For a similar setting, see, for example, \cite{RefC97}.
Let $\boldsymbol{Z}\in\mathbb{R}^{20\times17}$ satisfy $\boldsymbol{X}^{\top}\boldsymbol{Z} = \boldsymbol{0}$ and $\boldsymbol{Z}^{\top}\boldsymbol{Z} = \boldsymbol{I}_{17}$, and define
\begin{align*}
\boldsymbol{\Omega}_{c}(\tau) = \boldsymbol{I}_{20} + \tau \boldsymbol{Z}\boldsymbol{Z}^{\top} + c\tau \boldsymbol{P}_{X}, \quad c\in\{0,1\}.
\end{align*}
Note that the equality \eqref{B2EK} holds for any $\tau>0$ and $\boldsymbol{K} \in \MS^N(k)\setminus{\{\boldsymbol{0}\}}$ if and only if $c=0$.
We set 
\begin{align*}
\boldsymbol{\beta}
=
(2.0, 1.6, -0.8)^\top, \ \sigma^2 = 1, \ \boldsymbol{K}
= \diag(1.0,0.3,1.2),
\end{align*}
and consider $m\in\{30,50,100\}$ and $\tau\in\{1,4,8\}$.
For each combination of $m$, $\tau$, and $c$, the simulation procedure is as follows.
\begin{enumerate}
\item Generate $m$ independent training observations $\boldsymbol{y}_1, \ldots, \boldsymbol{y}_m$ from $N_{20}(\boldsymbol{X}\boldsymbol{\beta}, \boldsymbol{\Omega}_{c}(\tau))$, where $N_d(\boldsymbol{\mu},\boldsymbol{\Sigma})$ denotes the $d$-dimensional normal distribution with mean vector $\boldsymbol{\mu}\in\mathbb{R}^d$ and covariance matrix $\boldsymbol{\Sigma}\in\MS^+(d)$.

\item Estimate $\tau$ by the method of moments as
\begin{align*}
\hat{\tau}
=
\max\left\{
0,\,
\frac{1}{17m}
\sum_{j=1}^{m}
\left\|
\boldsymbol{Z}^{\top}\boldsymbol{y}_j
\right\|^2
-1
\right\}.
\end{align*}
and construct $\hat{\boldsymbol{\Omega}}_{c}
=
\boldsymbol{I}_{20}
+
\hat{\tau}
\boldsymbol{Z}\boldsymbol{Z}^{\top}
+
c\hat{\tau}
\boldsymbol{P}_{X}$.

\item Generate an independent test observation $\boldsymbol{y}_0$ from
$N_{20}(\boldsymbol{X}\boldsymbol{\beta}, \boldsymbol{\Omega}_{c}(\tau))$
and compute
\begin{align*}
\hat{\boldsymbol{\beta}}_{BL}
(\boldsymbol{y}_0;\hat{\boldsymbol{\Omega}}_{c})
=
\boldsymbol{L}
(\hat{\boldsymbol{\Omega}}_{c})
\boldsymbol{y}_0,
\qquad
\hat{\boldsymbol{\beta}}_{BL}
(\boldsymbol{y}_0;\boldsymbol{I}_{20})
=
\boldsymbol{L}
(\boldsymbol{I}_{20})
\boldsymbol{y}_0,
\end{align*}
where $\boldsymbol{L}(\boldsymbol{\Phi})
=
\boldsymbol{K}\boldsymbol{X}^{\top}
(
\boldsymbol{\Phi}
+
\boldsymbol{X}\boldsymbol{K}\boldsymbol{X}^{\top}
)^{-1}$.
\end{enumerate}
We repeat the procedure $R=5000$ times using a fixed random seed.
The estimation accuracy of $\hat{\tau}$ is evaluated by
\begin{align*}
\operatorname{RMSE}(\hat{\tau})
=
\left\{
\frac{1}{R}
\sum_{r=1}^{R}
\left(
\hat{\tau}^{(r)}-\tau
\right)^2
\right\}^{1/2}.
\end{align*}
The discrepancies between two BLEs are measured by
\begin{align*}
d_{L} = \frac{1}{R} \sum_{r=1}^{R} \left\|\boldsymbol{L}\left(\hat{\boldsymbol{\Omega}}_{c}^{(r)}\right) -
\boldsymbol{L}\left(\boldsymbol{I}_{20}\right)\right\|_\mathrm{F}, \quad
d_{\beta}
= \frac{1}{R}\sum_{r=1}^{R}\left\|\hat{\boldsymbol{\beta}}_{BL}^{(r)}(\boldsymbol{y}_0;\hat{\boldsymbol{\Omega}}_{c}) - \hat{\boldsymbol{\beta}}_{BL}^{(r)}(\boldsymbol{y}_0;\boldsymbol{I}_{20})
\right\|^2,
\end{align*}
where $\|\cdot\|_\mathrm{F}$ denotes the Frobenius norm.

\begin{table}[t]
\centering
\setlength{\belowcaptionskip}{8pt}
\caption{Monte Carlo results for the mixed-effects model}
\label{table2}
\small
\setlength{\tabcolsep}{2pt}
\begin{tabular}{cccccccc}
\toprule
& & \multicolumn{3}{c}{$c=0$} & \multicolumn{3}{c}{$c=1$} \\
\cmidrule(lr){3-5}\cmidrule(lr){6-8}
$m$
& $\tau$
& $\operatorname{RMSE}(\hat{\tau})$
& $d_L$
& $d_{\beta}$
& $\operatorname{RMSE}(\hat{\tau})$
& $d_L$
& $d_{\beta}$ \\
\midrule
30 & 1 & 0.1231 & $5.78\times10^{-16}$ & $7.68\times10^{-31}$ & 0.1257 & 0.1340 & 0.1762 \\
30 & 4 & 0.3083 & $5.83\times10^{-16}$ & $1.76\times10^{-30}$ & 0.3139 & 0.3132 & 1.2652 \\
30 & 8 & 0.5561 & $5.99\times10^{-16}$ & $3.18\times10^{-30}$ & 0.5630 & 0.4041 & 2.7571 \\
50 & 1 & 0.0975 & $5.79\times10^{-16}$ & $7.63\times10^{-31}$ & 0.0969 & 0.1340 & 0.1755 \\
50 & 4 & 0.2426 & $5.83\times10^{-16}$ & $1.76\times10^{-30}$ & 0.2426 & 0.3132 & 1.2804 \\
50 & 8 & 0.4434 & $6.00\times10^{-16}$ & $3.32\times10^{-30}$ & 0.4327 & 0.4043 & 2.7899 \\
100 & 1 & 0.0674 & $5.79\times10^{-16}$ & $7.79\times10^{-31}$ & 0.0680 & 0.1341 & 0.1775 \\
100 & 4 & 0.1748 & $5.83\times10^{-16}$ & $1.80\times10^{-30}$ & 0.1728 & 0.3136 & 1.2609 \\
100 & 8 & 0.3092 & $6.00\times10^{-16}$ & $3.32\times10^{-30}$ & 0.3069 & 0.4043 & 2.8130 \\
\bottomrule
\end{tabular}
\end{table}

As shown in Table~\ref{table2}, $\hat{\tau}$ has a non-negligible estimation error, especially when $m$ is small or $\tau$ is large.
Nevertheless, when $c=0$, both $d_L$ and $d_{\beta}$ remain at machine precision for all parameter settings.
Indeed, every realization of $\hat{\tau}$ satisfies
\begin{align*}
\hat{\boldsymbol{\Omega}}_{0}\boldsymbol{X}
= (
\boldsymbol{I}_{20}
+
\hat{\tau}
\boldsymbol{Z}\boldsymbol{Z}^{\top})
\boldsymbol{X} = \boldsymbol{X}.
\end{align*}
Consequently, estimation error in $\hat{\tau}$ does not affect the BLE, and
\begin{align*}
\hat{\boldsymbol{\beta}}_{BL}(\boldsymbol{y}_0;\hat{\boldsymbol{\Omega}}_0) \simeq \hat{\boldsymbol{\beta}}_{BL}(\boldsymbol{y}_0;\boldsymbol{I}_{20}).
\end{align*}
When $c=1$, both discrepancies are positive and increase with $\tau$, reflecting the increasing departure from the condition
$\boldsymbol{\Omega}_1(\tau)\boldsymbol{XK} = \boldsymbol{XK}$.
These results numerically support Theorem~\ref{thm2}.
\end{exam}
\newpage
\begin{exam}[\textbf{Spatial error model \citep{RefD82}}]\label{ex3}
\mbox{}\\
Consider the spatial error model
\begin{align}\label{ar1}
\left\{ \,
 \begin{aligned}
&\boldsymbol{y} = \boldsymbol{X}\boldsymbol{\beta} + \boldsymbol{u}\\ &\boldsymbol{u} = \rho \boldsymbol{W}\boldsymbol{u} + \boldsymbol{\ve}
\end{aligned}
    \right., \quad \E[\boldsymbol{\ve}] = \boldsymbol{0}_n, \quad \V(\boldsymbol{\ve}) = \sigma^2 \boldsymbol{I}_n,
\end{align}
where $\rho \in (-1,1)$ denotes an unknown spatial correlation coefficient for a given area partitioned into $n$ nonoverlapping regions
$R_i \ (i=1,\ldots,n)$, and $\boldsymbol{W} \in \mathbb{R}^{n \times n}$ is some known nonzero matrix of nonnegative weights. 
More specifically, we define a binary matrix $\boldsymbol{C} \in \mathbb{R}^{n \times n}$, whose elements are
\begin{align*}
c_{ij}=
\begin{cases}
1 & \textup{if regions $R_i$ and $R_j$ are contiguous}\\
0 & \textup{if regions $R_i$ and $R_j$ are not contiguous or $i=j$}
\end{cases},
\end{align*}
and $\boldsymbol{W}$ is formed using as weights $w_{ij} = c_{ij}/\sum_{j=1}^n c_{ij}$.
Thus, each row of $\boldsymbol{W}$ sums to 1, and the weights are simple proportions based on the number of neighboring regions for each area.
Such models are often used for data on plant heights, country populations, or positions in social networks \citep{RefKD87}.

The disturbance vector in \eqref{ar1} can be written as $\boldsymbol{u} = (\boldsymbol{I}_n - \rho \boldsymbol{W})^{-1}\boldsymbol{\ve}$.
Thus, we obtain the covariance matrix $\V(\boldsymbol{u}) = \sigma^2 \boldsymbol{\Omega}$ with $\boldsymbol{\Omega} = (\boldsymbol{I}_n - \rho \boldsymbol{W})^{-1}(\boldsymbol{I}_n - \rho \boldsymbol{W}^\top)^{-1}$.
\begin{cor}\label{cor410}
If $\MC(\boldsymbol{K}) \subseteq \mathcal{N}(\boldsymbol{W}\boldsymbol{X}) \cap \mathcal{N}(\boldsymbol{W}^\top \boldsymbol{X})$ in \eqref{ar1}, then the equality \eqref{B2EK} is satisfied.
\end{cor}
\begin{proof}
Since $\boldsymbol{W}\boldsymbol{X}\boldsymbol{K} = \boldsymbol{0}$ and $\boldsymbol{W}^\top\boldsymbol{X}\boldsymbol{K} = \boldsymbol{0}$, we have
\begin{align*}
\boldsymbol{\Omega}^{-1}\boldsymbol{XK} &=  (\boldsymbol{I}_n - \rho \boldsymbol{W}^\top)(\boldsymbol{I}_n - \rho \boldsymbol{W})\boldsymbol{XK}\\
& = \boldsymbol{X}\boldsymbol{K} - \rho (\boldsymbol{W} + \boldsymbol{W}^\top - \rho\boldsymbol{W}^\top\boldsymbol{W})\boldsymbol{X}\boldsymbol{K}\\
& = \boldsymbol{X}\boldsymbol{K}.
\end{align*}
Hence, it follows from \eqref{B2Enew} that the equality \eqref{B2EK} is satisfied.
\end{proof}

\begin{rem} 
(i) If the condition in Corollary~\ref{cor410} is satisfied, one can still determine whether the equality condition holds even if $\rho$ remains unknown.\\
(ii) \cite{RefG01} introduced a sufficient condition for the equivalence of $\hat{\boldsymbol{\beta}}_{GLS}$ and $\hat{\boldsymbol{\beta}}_{OLS}$ in a framework similar to that of Example~\ref{ex3}.
\end{rem}
\textbf{Numerical Experiment.}
We conduct a numerical experiment to examine the equivalence of two BLEs under \eqref{ar1}.
Let $n=20$ and $k=2$.
Suppose that the spatial units are arranged on a circle, and define the symmetric weight matrix $\boldsymbol{W}\in\mathbb{R}^{20\times20}$ by
\begin{align*}
w_{ij}
=
\begin{cases}
\dfrac{1}{2} & j\equiv i-1 \pmod{20}
       \textup{ or } j\equiv i+1 \pmod{20}\\
0 & \textup{otherwise}
\end{cases},
\quad
i,j=1,\ldots,20.
\end{align*}
For $i=1,\ldots,20$, let
\begin{align*}
x_{i1}
=
\cos\left\{\frac{\pi(i-1)}{2}\right\},
\qquad
x_{i2}
=
\sin\left\{\frac{\pi(i-1)}{2}\right\},
\end{align*}
and define $\boldsymbol{X}=(x_{ij})\in\mathbb{R}^{20\times2}$.
It holds that
\begin{align}\label{WXzero}
\boldsymbol{W}\boldsymbol{X}
=
\boldsymbol{W}^{\top}\boldsymbol{X}
=
\boldsymbol{0}.
\end{align}

The training observations $\boldsymbol{y}_1, \ldots, \boldsymbol{y}_m$ and an independent test observation $\boldsymbol{y}_0$ are generated from
\begin{align*}
\boldsymbol{y}_j
=
\boldsymbol{X}\boldsymbol{\beta}
+
\left(
\boldsymbol{I}_{20}
-
\rho\boldsymbol{W}
\right)^{-1}
\boldsymbol{\ve}_j,
\qquad
j=0,1,\ldots,m,
\end{align*}
where
$\boldsymbol{\ve}_j\sim
N_{20}(\boldsymbol{0}_{20},\sigma^2\boldsymbol{I}_{20})$,
and all random components are mutually independent.
Thus, the covariance matrix of $\boldsymbol{y}_j$ is
$\sigma^2\boldsymbol{\Omega}(\rho)$, where $\boldsymbol{\Omega}(\rho)
=
(
\boldsymbol{I}_{20}
-
\rho\boldsymbol{W}
)^{-1}
(
\boldsymbol{I}_{20}
-
\rho\boldsymbol{W}^{\top}
)^{-1}$.
We set
\begin{align*}
\boldsymbol{\beta}=(1.0,0.8)^{\top}, \ \sigma^2=1, \ \boldsymbol{K}=\operatorname{diag}(1.0,0.7),
\end{align*}
and consider $m\in\{30,50,100\}$ and $\rho\in\{0.2,0.5,0.8\}$.

For each combination of $m$ and $\rho$, we follow the simulation procedure in Example~\ref{ex2}, with $\rho$ estimated by maximum likelihood. 
For a candidate value $\rho \in[0,1)$, define
\begin{align*}
\hat{\boldsymbol{\beta}}(\rho)
=
(\boldsymbol{X}^{\top}
\boldsymbol{\Omega}(\rho)^{-1}
\boldsymbol{X})^{-1}
\boldsymbol{X}^{\top}
\boldsymbol{\Omega}(\rho)^{-1}
\overline{\boldsymbol{y}},
\end{align*}
where 
\begin{align*}
\overline{\boldsymbol{y}}
=
\frac{1}{m}
\sum_{j=1}^{m}
\boldsymbol{y}_j \quad \an \quad 
\hat{\sigma}^2(\rho)
=
\frac{1}{20m}
\sum_{j=1}^{m}
\left(
\boldsymbol{y}_j
-
\boldsymbol{X}\hat{\boldsymbol{\beta}}(\rho)
\right)^{\top}
\boldsymbol{\Omega}(\rho)^{-1}
\left(
\boldsymbol{y}_j
-
\boldsymbol{X}\hat{\boldsymbol{\beta}}(\rho)
\right).
\end{align*}
The maximum likelihood estimator of $\rho$ is given by
\begin{align*}
\hat{\rho}
=
\operatorname*{argmin}_{0\leq \rho < 1}
\left[
\log
\left|
\boldsymbol{\Omega}(\rho)
\right|
+
20\log
\left(
\hat{\sigma}^2(\rho)
\right)
\right],
\end{align*}
and we construct
$\hat{\boldsymbol{\Omega}}=\boldsymbol{\Omega}(\hat{\rho})$.

\begin{table}[t]
\centering
\setlength{\belowcaptionskip}{8pt}
\caption{Monte Carlo results for the spatial error model}
\label{table3}
\begin{tabular}{ccccc}
\toprule
$m$
& $\rho$
& $\operatorname{RMSE}(\hat{\rho})$
& $d_L$
& $d_{\beta}$\\
\midrule
30 & 0.2 & 0.0397 & $4.93\times10^{-16}$ & $7.11\times10^{-31}$\\
30 & 0.5 & 0.0319 & $6.96\times10^{-16}$ & $1.31\times10^{-30}$\\
30 & 0.8 & 0.0171 & $1.26\times10^{-15}$ & $4.56\times10^{-30}$\\
50 & 0.2 & 0.0306 & $4.87\times10^{-16}$ & $6.80\times10^{-31}$\\
50 & 0.5 & 0.0248 & $6.93\times10^{-16}$ & $1.34\times10^{-30}$\\
50 & 0.8 & 0.0133 & $1.26\times10^{-15}$ & $4.54\times10^{-30}$\\
100 & 0.2 & 0.0214 & $4.83\times10^{-16}$ & $6.64\times10^{-31}$\\
100 & 0.5 & 0.0174 & $6.91\times10^{-16}$ & $1.39\times10^{-30}$\\
100 & 0.8 & 0.0093 & $1.26\times10^{-15}$ & $4.23\times10^{-30}$\\
\bottomrule
\end{tabular}
\end{table}

Table~\ref{table3} shows that RMSE of $\hat{\rho}$ decreases as $m$ increases.
However, both $d_L$ and $d_{\beta}$ remain at machine precision for all settings.
This is because \eqref{WXzero} implies that, for every value of $\hat{\rho}$,
\begin{align*}
\hat{\boldsymbol{\Omega}}^{-1}
\boldsymbol{X}\boldsymbol{K}
= (\boldsymbol{I}_{20} - \hat{\rho}\boldsymbol{W}^{\top})
(\boldsymbol{I}_{20} - \hat{\rho}\boldsymbol{W}) \boldsymbol{X}\boldsymbol{K} =
\boldsymbol{X}\boldsymbol{K}.
\end{align*}
Thus, $\hat{\boldsymbol{\beta}}_{BL}(\boldsymbol{I}_{20},\boldsymbol{K})$ yields the same estimate as $\hat{\boldsymbol{\beta}}_{BL}(\hat{\boldsymbol{\Omega}}, \boldsymbol{K})$ without estimating the unknown spatial correlation coefficient.
\end{exam}

\section{Set of \texorpdfstring{$\boldsymbol{y}$}{y} satisfying the equality between two Bayes linear estimators}\label{sec:5}
In this section, we characterize the set of $\boldsymbol{y}$ for which the equality \eqref{B2Ey} holds.
The following lemma will be used in the proof of Theorem~\ref{thm3}.
\begin{lem}\label{lem43}
Let $\boldsymbol{A} \in \mathbb{R}^{q \times n}$ and $\boldsymbol{Q}_X = \boldsymbol{I}_n - \boldsymbol{P}_X$.
If $\boldsymbol{A}\boldsymbol{Q}_X = \boldsymbol{0}$, then $\mathcal{N}(\boldsymbol{A}) = \{\MC(\boldsymbol{X}) \cap \mathcal{N}(\boldsymbol{A})\} \oplus \MC(\boldsymbol{X})^\bot$,
where $\mathcal{N}(\boldsymbol{M})$ denotes the null space of a matrix $\boldsymbol{M}$.
\end{lem}
\begin{proof}
For all $\boldsymbol{y} \in \mathbb{R}^n$, $\boldsymbol{y}$ can be uniquely decomposed as $\boldsymbol{y} = \boldsymbol{P}_X \boldsymbol{y} + \boldsymbol{Q}_X\boldsymbol{y}$, where $\boldsymbol{P}_X\boldsymbol{y} \in \MC(\boldsymbol{X})$ and $\boldsymbol{Q}_X\boldsymbol{y} \in \MC(\boldsymbol{X})^\bot$.
Since $\boldsymbol{A}\boldsymbol{Q}_X=\boldsymbol{0}$, we have $\boldsymbol{A}\boldsymbol{y} = \boldsymbol{A}\boldsymbol{P}_X\boldsymbol{y}$.
Hence, $\boldsymbol{y} \in \mathcal{N}(\boldsymbol{A})$ is equivalent to
$\boldsymbol{P}_X\boldsymbol{y} \in \MC(\boldsymbol{X}) \ \cap \ \mathcal{N}(\boldsymbol{A})$.
Therefore, $\mathcal{N}(\boldsymbol{A}) = \{\MC(\boldsymbol{X}) \ \cap \ \mathcal{N}(\boldsymbol{A})\} \oplus \MC(\boldsymbol{X})^\bot$, where the directness follows from $\MC(\boldsymbol{X}) \cap \MC(\boldsymbol{X})^\bot = \{0\}$.
This completes the proof.
\end{proof}

Under \eqref{Rcond}, the set of $\boldsymbol{y}$ satisfying the equality \eqref{B2Ey} can be written as a direct sum decomposition.
The following theorem is one of the main results of this paper.
\begin{thm}\label{thm3}
Suppose that $\boldsymbol{K}_1, \boldsymbol{K}_2 \in  \MS^N(k)$ are fixed and 
$\boldsymbol{\Omega}$ is of the form \eqref{Rcond}.
Then the equality \eqref{B2Ey} holds if and only if
\begin{align}\label{oplus}
\boldsymbol{y} \in \left\{\mathcal{C}(\boldsymbol{X}) \cap \mathcal{N}\left(\boldsymbol{A}\right)\right\} \oplus \MC(\boldsymbol{X})^\bot,
\end{align}
where $\boldsymbol{A} = (\boldsymbol{K}_2 \boldsymbol{X}^\top \boldsymbol{\Omega} - \boldsymbol{K}_1 \boldsymbol{X}^\top)(\boldsymbol{\Omega}+ \boldsymbol{X}\boldsymbol{K}_1\boldsymbol{X}^\top)^{-1}$.
\end{thm}
\begin{proof}
To see the necessity, suppose that \eqref{B2Ey} is satisfied.
Then \eqref{B2Ey} is equivalent to
\begin{align*}
&\hat{\boldsymbol{\beta}}_{BL}(\boldsymbol{\Omega}, \boldsymbol{K}_1) = \hat{\boldsymbol{\beta}}_{BL} (\boldsymbol{I}_n, \boldsymbol{K}_2) \\
\Leftrightarrow \quad & \boldsymbol{\Omega}(\boldsymbol{\Omega}+\boldsymbol{X}\boldsymbol{K}_1\boldsymbol{X}^\top)^{-1}\boldsymbol{y} = (\boldsymbol{I}_n+\boldsymbol{X}\boldsymbol{K}_2\boldsymbol{X}^\top)^{-1}\boldsymbol{y}\\
\Leftrightarrow \quad & \{(\boldsymbol{I}_n+\boldsymbol{X}\boldsymbol{K}_2\boldsymbol{X}^\top)\boldsymbol{\Omega}(\boldsymbol{\Omega}+\boldsymbol{X}\boldsymbol{K}_1\boldsymbol{X}^\top)^{-1} - \boldsymbol{I}_n\}\boldsymbol{y} = \boldsymbol{0}_n \notag\\
\Leftrightarrow \quad & \{(\boldsymbol{I}_n+\boldsymbol{X}\boldsymbol{K}_2\boldsymbol{X}^\top)\boldsymbol{\Omega} - (\boldsymbol{\Omega}+\boldsymbol{X}\boldsymbol{K}_1\boldsymbol{X}^\top)\}(\boldsymbol{\Omega}+\boldsymbol{X}\boldsymbol{K}_1\boldsymbol{X}^\top)^{-1}\boldsymbol{y} = \boldsymbol{0}_n \\
\Leftrightarrow \quad & \boldsymbol{X}(\boldsymbol{K}_2\boldsymbol{X}^\top \boldsymbol{\Omega} - \boldsymbol{K}_1\boldsymbol{X}^\top)(\boldsymbol{\Omega}+\boldsymbol{X}\boldsymbol{K}_1\boldsymbol{X}^\top)^{-1}\boldsymbol{y} = \boldsymbol{0}_n \\
\Leftrightarrow \quad & \boldsymbol{A}\boldsymbol{y} = \boldsymbol{0}_k\\
\Leftrightarrow \quad & \boldsymbol{y} \in \mathcal{N}(\boldsymbol{A}).
\end{align*}
It follows from \eqref{Rcond} that $\boldsymbol{\Omega}\boldsymbol{Q}_X = \boldsymbol{Q}_X \boldsymbol{\Omega}$.
Additionally, it holds that $\boldsymbol{Q}_X\boldsymbol{X}\boldsymbol{K}_1\boldsymbol{X}^\top = \boldsymbol{X}\boldsymbol{K}_1\boldsymbol{X}^\top\boldsymbol{Q}_X$.
So, we have
\begin{align*}
&(\boldsymbol{\Omega}+\boldsymbol{X}\boldsymbol{K}_1\boldsymbol{X}^\top) \boldsymbol{Q}_X = \boldsymbol{Q}_X (\boldsymbol{\Omega}+\boldsymbol{X}\boldsymbol{K}_1\boldsymbol{X}^\top)\\ 
\Leftrightarrow \quad &  (\boldsymbol{\Omega}+\boldsymbol{X}\boldsymbol{K}_1\boldsymbol{X}^\top)^{-1}\boldsymbol{Q}_X = \boldsymbol{Q}_X(\boldsymbol{\Omega}+\boldsymbol{X}\boldsymbol{K}_1\boldsymbol{X}^\top)^{-1}.
\end{align*}
Thus, we obtain
\begin{align*}
\boldsymbol{A}\boldsymbol{Q}_X & = (\boldsymbol{K}_2 \boldsymbol{X}^\top \boldsymbol{\Omega} - \boldsymbol{K}_1 \boldsymbol{X}^\top)(\boldsymbol{\Omega}+ \boldsymbol{X}\boldsymbol{K}_1\boldsymbol{X}^\top)^{-1}\boldsymbol{Q}_X\\
& = (\boldsymbol{K}_2 \boldsymbol{X}^\top \boldsymbol{\Omega} - \boldsymbol{K}_1 \boldsymbol{X}^\top) \boldsymbol{Q}_X (\boldsymbol{\Omega}+ \boldsymbol{X}\boldsymbol{K}_1\boldsymbol{X}^\top)^{-1}\\
& = (\boldsymbol{K}_2 \boldsymbol{X}^\top \boldsymbol{\Omega}\boldsymbol{Q}_X - \boldsymbol{K}_1 \boldsymbol{X}^\top\boldsymbol{Q}_X)  (\boldsymbol{\Omega}+ \boldsymbol{X}\boldsymbol{K}_1\boldsymbol{X}^\top)^{-1}\\
& = (\boldsymbol{K}_2 \boldsymbol{X}^\top\boldsymbol{Q}_X \boldsymbol{\Omega} - \boldsymbol{K}_1 \boldsymbol{X}^\top\boldsymbol{Q}_X)  (\boldsymbol{\Omega}+ \boldsymbol{X}\boldsymbol{K}_1\boldsymbol{X}^\top)^{-1}\\
& = \boldsymbol{0}.
\end{align*}
From Lemma~\ref{lem43}, the equality \eqref{B2Ey} is equivalent to \eqref{oplus}.
This completes the proof.
\end{proof}

\begin{rem}
(i) Suppose that $\boldsymbol{K}_1, \boldsymbol{K}_2 \in  \MS^N(k)$ are fixed and 
$\boldsymbol{\Omega}$ is of the form \eqref{Rcond}.
Then the equality \eqref{B2E} holds if and only if
\begin{align*}
\mathcal{C}(\boldsymbol{X}) \subseteq \mathcal{N}\left[(\boldsymbol{K}_2 \boldsymbol{X}^\top \boldsymbol{\Omega} - \boldsymbol{K}_1 \boldsymbol{X}^\top)(\boldsymbol{\Omega}+ \boldsymbol{X}\boldsymbol{K}_1\boldsymbol{X}^\top)^{-1}\right].
\end{align*}
This proof can be carried out by using the decomposition $\mathbb{R}^n = \mathcal{C}(\boldsymbol{X})\oplus\MC(\boldsymbol{X})^\bot$.\\
(ii) \cite{RefG04} showed that $\MC(\begin{matrix}
\boldsymbol{\Omega} &  \boldsymbol{X}
\end{matrix}) = \MC(\boldsymbol{X})\oplus \mathcal{C}(\boldsymbol{\Omega}\boldsymbol{Z})$. 
This implies that the necessary and sufficient condition \eqref{oplus} can also be written as 
\begin{align*}
\boldsymbol{y} \in \left\{\mathcal{C}(\boldsymbol{X})\ \cap \ \mathcal{N}\left[(\boldsymbol{K}_2 \boldsymbol{X}^\top \boldsymbol{\Omega} - \boldsymbol{K}_1 \boldsymbol{X}^\top)(\boldsymbol{\Omega}+ \boldsymbol{X}\boldsymbol{K}_1\boldsymbol{X}^\top)^{-1}\right]\right\}  \oplus \mathcal{C}(\boldsymbol{\Omega Z}).   
\end{align*}
(iii) Let us consider Rao's mixed-effects model; see Example~\ref{ex2}.
Then $\boldsymbol{\Omega}$ is of the form \eqref{Rcond}, and thus the assumption of Theorem~\ref{thm3} is satisfied.\\
(iv) The necessary and sufficient condition for the equality \eqref{B2Ey} is
\begin{align}\label{yy}
\boldsymbol{y} \in \mathcal{N}\left[(\boldsymbol{K}_2 \boldsymbol{X}^\top \boldsymbol{\Omega} - \boldsymbol{K}_1 \boldsymbol{X}^\top)(\boldsymbol{\Omega}+ \boldsymbol{X}\boldsymbol{K}_1\boldsymbol{X}^\top)^{-1}\right]
\end{align}
for given $\boldsymbol{K}_1,\boldsymbol{K}_2 \in \mathcal{S}^{N}(k)$, and we can directly derive Theorem~\ref{thm2} from \eqref{yy}.
In fact, it holds that
\allowdisplaybreaks\begin{align*}
&\hat{\boldsymbol{\beta}}_{BL}(\boldsymbol{\Omega},\boldsymbol{K}_1) \equiv \hat{\boldsymbol{\beta}}_{BL}(\boldsymbol{I}_n,\boldsymbol{K}_2)\\
\Leftrightarrow \quad & (\boldsymbol{K}_2 \boldsymbol{X}^\top \boldsymbol{\Omega} - \boldsymbol{K}_1 \boldsymbol{X}^\top)(\boldsymbol{\Omega}+\boldsymbol{X}\boldsymbol{K}_1\boldsymbol{X}^\top)^{-1}\boldsymbol{y} \equiv \boldsymbol{0}_k\\
\Leftrightarrow \quad & (\boldsymbol{K}_2 \boldsymbol{X}^\top \boldsymbol{\Omega} - \boldsymbol{K}_1 \boldsymbol{X}^\top)(\boldsymbol{\Omega}+\boldsymbol{X}\boldsymbol{K}_1\boldsymbol{X}^\top)^{-1} =\boldsymbol{0}\\
\Leftrightarrow \quad & \boldsymbol{K}_2 \boldsymbol{X}^\top \boldsymbol{\Omega} = \boldsymbol{K}_1 \boldsymbol{X}^\top,
\end{align*}
and we have \eqref{4T}.\\
(v) Suppose that $\boldsymbol{K} \in \MS^+(k)$ is fixed and $\boldsymbol{\Omega}$ is of the form \eqref{Rcond}.
Then the equality $\hat{\boldsymbol{\beta}}_{BL}(\boldsymbol{\Omega}, \boldsymbol{K}) = \hat{\boldsymbol{\beta}}_{BL} (\boldsymbol{I}_n, \boldsymbol{K})$ holds if and only if 
\begin{align*}
\boldsymbol{y} \in \left\{\MC(\boldsymbol{X})\ \cap \ \mathcal{N}\left[\boldsymbol{X}^\top(\boldsymbol{\Omega} - \boldsymbol{I}_n)(\boldsymbol{\Omega}+ \boldsymbol{X}\boldsymbol{K}\boldsymbol{X}^\top)^{-1}\right]\right\}  \oplus \MC(\boldsymbol{X})^\bot.
\end{align*}
\end{rem}

\section{Equivalence of residual sums of squares under Bayes linear estimators}\label{sec:6}
In this section, we mainly consider the necessary and sufficient conditions under which the two residual sums of squares coincide for all $\boldsymbol{y} \in \mathbb{R}^n$.

To simplify the notation, let us denote
\begin{align*}
\boldsymbol{S} = (\boldsymbol{\Omega}+\boldsymbol{X}\boldsymbol{K}_1\boldsymbol{X}^{\top})^{-1}, \quad \boldsymbol{T}= (\boldsymbol{I}_n+\boldsymbol{X}\boldsymbol{K}_2\boldsymbol{X}^{\top})^{-1}.
\end{align*}
The following theorem is one of the main results of this paper.
\begin{thm}\label{thm4}
Let $\boldsymbol{\Omega}$ be written as \eqref{Orep}.
For given $\boldsymbol{K}_1, \boldsymbol{K}_2 \in \MS^N(k)$, the equality \eqref{BRSS1} holds
if and only if the two conditions 
\begin{align}
& \boldsymbol{\Omega}= \boldsymbol{X}\boldsymbol{\Gamma} \boldsymbol{X}^{\top}+\boldsymbol{Z}(\boldsymbol{Z}^{\top}\boldsymbol{Z})^{-1}\boldsymbol{Z}^{\top}, \label{516}\\
& (\boldsymbol{\Gamma} +\boldsymbol{K}_1)\boldsymbol{\Gamma}^{-1} (\boldsymbol{\Gamma} +\boldsymbol{K}_1)=\left[(\boldsymbol{X}^{\top}\boldsymbol{X})^{-1}+\boldsymbol{K}_2\right]\boldsymbol{X}^{\top}\boldsymbol{X}\left[(\boldsymbol{X}^{\top}\boldsymbol{X})^{-1}+\boldsymbol{K}_2\right]\label{517}
\end{align}
simultaneously hold for some $\boldsymbol{\Gamma} \in \mathcal{S}^+(k)$.
\end{thm}
\begin{proof}
Since $\boldsymbol{y}-\boldsymbol{X}\hat{\boldsymbol{\beta}}_{BL}(\boldsymbol{\Omega},\boldsymbol{K}_1) 
=\boldsymbol{\Omega} \boldsymbol{S}\boldsymbol{y}$ and $\boldsymbol{y}-\boldsymbol{X}\hat{\boldsymbol{\beta}}_{BL}(\boldsymbol{I}_n,\boldsymbol{K}_2) =\boldsymbol{T}\boldsymbol{y},$
the residual sums of squares $RSS_{BL}(\boldsymbol{\Omega},\boldsymbol{K}_1)$ and $RSS_{BL}(\boldsymbol{I}_n,\boldsymbol{K}_2)$ can be rewritten as
\begin{align*}
&RSS_{BL}(\boldsymbol{\Omega},\boldsymbol{K}_1)=\left(\boldsymbol{y}-\boldsymbol{X}\hat{\boldsymbol{\beta}}_{BL}(\boldsymbol{\Omega},\boldsymbol{K}_1)\right)^{\top}\boldsymbol{\Omega}^{-1}\left(\boldsymbol{y}-\boldsymbol{X}\hat{\boldsymbol{\beta}}_{BL}(\boldsymbol{\Omega},\boldsymbol{K}_1)\right)=\boldsymbol{y}^{\top}\boldsymbol{S}\boldsymbol{\Omega}\boldsymbol{S}\boldsymbol{y}
\end{align*}
and 
\begin{align*}
RSS_{BL}(\boldsymbol{I}_n,\boldsymbol{K}_2)=\left(\boldsymbol{y}-\boldsymbol{X}\hat{\boldsymbol{\beta}}_{BL}(\boldsymbol{I}_n,\boldsymbol{K}_2)\right)^{\top} \left(\boldsymbol{y}-\boldsymbol{X}\hat{\boldsymbol{\beta}}_{BL}(\boldsymbol{I}_n,\boldsymbol{K}_2)\right)=\boldsymbol{y}^{\top}\boldsymbol{T}\boldsymbol{T} \boldsymbol{y}.
\end{align*}
Thus, the problem is to derive a condition under which $\boldsymbol{y}^\top (\boldsymbol{S}\boldsymbol{\Omega} \boldsymbol{S}-\boldsymbol{T}\boldsymbol{T})\boldsymbol{y} \equiv 0$ holds, which is equivalent to deriving a condition for
\begin{align}\label{RSSBL}
\boldsymbol{S}^{-1}\boldsymbol{\Omega}^{-1} \boldsymbol{S}^{-1} -\boldsymbol{T}^{-1}\boldsymbol{T}^{-1}=\boldsymbol{0}.
\end{align}
The equality \eqref{RSSBL} holds if and only if the following three
equalities simultaneously hold:
\begin{align}
\boldsymbol{X}^\top (\boldsymbol{S}^{-1}\boldsymbol{\Omega}^{-1} \boldsymbol{S}^{-1}-\boldsymbol{T}^{-1}\boldsymbol{T}^{-1}  ) \boldsymbol{X} =\boldsymbol{0}, \label{XXK3} \\
\boldsymbol{X}^\top (\boldsymbol{S}^{-1}\boldsymbol{\Omega}^{-1} \boldsymbol{S}^{-1}-\boldsymbol{T}^{-1}\boldsymbol{T}^{-1}  )  \boldsymbol{Z} =\boldsymbol{0}, \label{XZK3} \\
\boldsymbol{Z}^\top (\boldsymbol{S}^{-1}\boldsymbol{\Omega}^{-1} \boldsymbol{S}^{-1}-\boldsymbol{T}^{-1}\boldsymbol{T}^{-1}  ) \boldsymbol{Z} =\boldsymbol{0}. \label{ZZK3}
\end{align}
The quantities  $\boldsymbol{X}^{\top}\boldsymbol{S}^{-1}\boldsymbol{\Omega}^{-1} \boldsymbol{S}^{-1}\boldsymbol{X}, \boldsymbol{X}^{\top}\boldsymbol{S}^{-1}\boldsymbol{\Omega}^{-1} \boldsymbol{S}^{-1}\boldsymbol{Z},$ and $\boldsymbol{Z}^{\top}\boldsymbol{S}^{-1}\boldsymbol{\Omega}^{-1} \boldsymbol{S}^{-1}\boldsymbol{Z}$ are calculated as
\begin{align*}
\boldsymbol{X}^{\top}\boldsymbol{S}^{-1}
\boldsymbol{\Omega}^{-1}
\boldsymbol{S}^{-1}\boldsymbol{X}
&=
\boldsymbol{X}^{\top}
(\boldsymbol{\Omega}+\boldsymbol{X}\boldsymbol{K}_1\boldsymbol{X}^{\top})
\boldsymbol{\Omega}^{-1}
(\boldsymbol{\Omega}+\boldsymbol{X}\boldsymbol{K}_1\boldsymbol{X}^{\top})
\boldsymbol{X} \\
&=
\boldsymbol{X}^{\top}\boldsymbol{\Omega}\boldsymbol{X}
+
2\boldsymbol{X}^{\top}\boldsymbol{X}\boldsymbol{K}_1\boldsymbol{X}^{\top}\boldsymbol{X}
+
\boldsymbol{X}^{\top}\boldsymbol{X}
\boldsymbol{K}_1
\boldsymbol{X}^{\top}\boldsymbol{\Omega}^{-1}\boldsymbol{X}
\boldsymbol{K}_1
\boldsymbol{X}^{\top}\boldsymbol{X} \\
&=
\boldsymbol{X}^{\top}\boldsymbol{X}
\boldsymbol{\Gamma}
\boldsymbol{X}^{\top}\boldsymbol{X}
+
2\boldsymbol{X}^{\top}\boldsymbol{X}
\boldsymbol{K}_1
\boldsymbol{X}^{\top}\boldsymbol{X}
+
\boldsymbol{X}^{\top}\boldsymbol{X}
\boldsymbol{K}_1
\boldsymbol{A}
\boldsymbol{K}_1
\boldsymbol{X}^{\top}\boldsymbol{X},
\\[1ex]
\boldsymbol{X}^{\top}\boldsymbol{S}^{-1}
\boldsymbol{\Omega}^{-1}
\boldsymbol{S}^{-1}\boldsymbol{Z}
&=
\boldsymbol{X}^{\top}
(\boldsymbol{\Omega}+\boldsymbol{X}\boldsymbol{K}_1\boldsymbol{X}^{\top})
\boldsymbol{\Omega}^{-1}
(\boldsymbol{\Omega}+\boldsymbol{X}\boldsymbol{K}_1\boldsymbol{X}^{\top})
\boldsymbol{Z} \\
&=
\boldsymbol{X}^{\top}\boldsymbol{\Omega}\boldsymbol{Z}
+
2\boldsymbol{X}^{\top}\boldsymbol{X}
\boldsymbol{K}_1
\boldsymbol{X}^{\top}\boldsymbol{Z}
+
\boldsymbol{X}^{\top}\boldsymbol{X}
\boldsymbol{K}_1
\boldsymbol{X}^{\top}\boldsymbol{\Omega}^{-1}\boldsymbol{X}
\boldsymbol{K}_1
\boldsymbol{X}^{\top}\boldsymbol{Z} \\
&=
\boldsymbol{X}^{\top}\boldsymbol{X}
\boldsymbol{\Xi}
\boldsymbol{Z}^{\top}\boldsymbol{Z},
\\[1ex]
\boldsymbol{Z}^{\top}\boldsymbol{S}^{-1}
\boldsymbol{\Omega}^{-1}
\boldsymbol{S}^{-1}\boldsymbol{Z}
&=
\boldsymbol{Z}^{\top}
(\boldsymbol{\Omega}+\boldsymbol{X}\boldsymbol{K}_1\boldsymbol{X}^{\top})
\boldsymbol{\Omega}^{-1}
(\boldsymbol{\Omega}+\boldsymbol{X}\boldsymbol{K}_1\boldsymbol{X}^{\top})
\boldsymbol{Z} \\
&=
\boldsymbol{Z}^{\top}\boldsymbol{\Omega}\boldsymbol{Z}
+
2\boldsymbol{Z}^{\top}\boldsymbol{X}
\boldsymbol{K}_1
\boldsymbol{X}^{\top}\boldsymbol{Z}
+
\boldsymbol{Z}^{\top}\boldsymbol{X}
\boldsymbol{K}_1
\boldsymbol{X}^{\top}\boldsymbol{\Omega}^{-1}\boldsymbol{X}
\boldsymbol{K}_1
\boldsymbol{X}^{\top}\boldsymbol{Z} \\
&=
\boldsymbol{Z}^{\top}\boldsymbol{Z}
\boldsymbol{\Delta}
\boldsymbol{Z}^{\top}\boldsymbol{Z},
\end{align*}
where $\boldsymbol{A} = (\boldsymbol{\Gamma} - \boldsymbol{\Xi}\boldsymbol{\Delta}^{-1}\boldsymbol{\Xi}^\top)^{-1}$.
A straightforward calculation yields that
\begin{align*}
\boldsymbol{T}^{-1}\boldsymbol{T}^{-1} 
=(\boldsymbol{I}_n+\boldsymbol{X}\boldsymbol{K}_2\boldsymbol{X}^{\top})(\boldsymbol{I}_n+\boldsymbol{X}\boldsymbol{K}_2\boldsymbol{X}^{\top})
=\boldsymbol{I}_n + 2\boldsymbol{X}\boldsymbol{K}_2\boldsymbol{X}^\top + \boldsymbol{X}\boldsymbol{K}_2\boldsymbol{X}^\top\boldsymbol{X}\boldsymbol{K}_2\boldsymbol{X}^\top.
\end{align*}
Noting that $\boldsymbol{X}^\top\boldsymbol{T}^{-1}\boldsymbol{T}^{-1}\boldsymbol{Z} = \boldsymbol{0}$ and $\boldsymbol{Z}^\top \boldsymbol{T}^{-1}\boldsymbol{T}^{-1}\boldsymbol{Z} = \boldsymbol{Z}^\top \boldsymbol{Z}$, \eqref{XZK3} can be rewritten as 
\begin{align}\label{XZK4}
&\boldsymbol{X}^\top \boldsymbol{X} \boldsymbol{\Xi}\boldsymbol{Z}^\top \boldsymbol{Z} =\boldsymbol{0} \notag \\
\Leftrightarrow \quad &  \boldsymbol{\Xi} =\boldsymbol{0},
\end{align}
and \eqref{ZZK3} is equivalent to
\begin{align}\label{ZZK4}
&\boldsymbol{Z}^\top \boldsymbol{Z} \boldsymbol{\Delta}\boldsymbol{Z}^\top \boldsymbol{Z} =\boldsymbol{Z}^\top \boldsymbol{Z} \notag\\
\Leftrightarrow \quad &\boldsymbol{\Delta} =(\boldsymbol{Z}^\top \boldsymbol{Z})^{-1}.
\end{align}
From \eqref{XZK4}, $\boldsymbol{A}$ becomes $\boldsymbol{\Gamma}^{-1}$. 
Hence, the equality \eqref{XXK3} is equivalent to
\begin{align*}
&\boldsymbol{X}^{\top}\boldsymbol{X}\boldsymbol{\Gamma} \boldsymbol{X}^{\top}\boldsymbol{X}  + 2\boldsymbol{X}^{\top}\boldsymbol{X}\boldsymbol{K}_1\boldsymbol{X}^{\top}\boldsymbol{X}+\boldsymbol{X}^{\top}\boldsymbol{X}\boldsymbol{K}_1 \boldsymbol{\Gamma}^{-1}\boldsymbol{K}_1\boldsymbol{X}^{\top}\boldsymbol{X} \\
&= \boldsymbol{X}^\top \boldsymbol{X} + 2\boldsymbol{X}^\top \boldsymbol{X}\boldsymbol{K}_2\boldsymbol{X}^\top \boldsymbol{X} + \boldsymbol{X}^\top \boldsymbol{X}\boldsymbol{K}_2\boldsymbol{X}^\top\boldsymbol{X}\boldsymbol{K}_2\boldsymbol{X}^\top \boldsymbol{X}\\
\Leftrightarrow \quad & \boldsymbol{\Gamma} + 2\boldsymbol{K}_1 + \boldsymbol{K}_1\boldsymbol{\Gamma}^{-1}\boldsymbol{K}_1 = (\boldsymbol{X}^\top \boldsymbol{X})^{-1} + 2\boldsymbol{K}_2 + \boldsymbol{K}_2\boldsymbol{X}^\top \boldsymbol{X}\boldsymbol{K}_2\\
\Leftrightarrow \quad & (\boldsymbol{\Gamma} +\boldsymbol{K}_1)\boldsymbol{\Gamma} ^{-1} (\boldsymbol{\Gamma} +\boldsymbol{K}_1)=\left[(\boldsymbol{X}^{\top}\boldsymbol{X})^{-1}+\boldsymbol{K}_2\right]\boldsymbol{X}^{\top}\boldsymbol{X}\left[(\boldsymbol{X}^{\top}\boldsymbol{X})^{-1}+\boldsymbol{K}_2\right].
\end{align*}
It follows from  \eqref{XZK4} and \eqref{ZZK4} that $\boldsymbol{\Omega}$ is of the form \eqref{516}.
This completes the proof.
\end{proof}

\begin{rem}
Let us denote the residual sums of squares of $\hat{\boldsymbol{\beta}}_{OLS}$ and $\hat{\boldsymbol{\beta}}_{GLS}$ by
\begin{align*}
RSS(\hat{\boldsymbol{\beta}}_{OLS}) = \left\|\boldsymbol{y} - \boldsymbol{X}\hat{\boldsymbol{\beta}}_{OLS}\right\|^2, \quad  RSS(\hat{\boldsymbol{\beta}}_{GLS}) = \left\| \boldsymbol{\Omega}^{-1/2} \left(\boldsymbol{y} - \boldsymbol{X}\hat{\boldsymbol{\beta}}_{GLS}\right) \right\|^2,
\end{align*}
respectively.
Theorem~\ref{thm4} implies that if the equality \eqref{BRSS1} holds, then the two equalities $\hat{\boldsymbol{\beta}}_{GLS} \equiv \hat{\boldsymbol{\beta}}_{OLS}$ and $RSS(\hat{\boldsymbol{\beta}}_{GLS}) \equiv RSS(\hat{\boldsymbol{\beta}}_{OLS})$ simultaneously hold.
See also \cite{RefK80}.
\end{rem}

Hereafter, we consider the case $\boldsymbol{K}_1,\boldsymbol{K}_2 \in \MS^+(k)$ and proceed to the discussion of Theorem~\ref{thm4}.

\begin{rem}[\textbf{Equality between two residual sums of squares of typical shrinkage estimators}]
Let $\boldsymbol{K}_1 = \rho(\boldsymbol{X}^\top \boldsymbol{\Omega}^{-1}\boldsymbol{X})^{-1}$ and $\boldsymbol{K}_2 = \rho(\boldsymbol{X}^\top \boldsymbol{X})^{-1}$, where $\rho$ is a positive constant.
Then the equality $RSS_{BL}(\boldsymbol{\Omega}, \rho(\boldsymbol{X}^\top \boldsymbol{\Omega}^{-1}\boldsymbol{X})^{-1}) \equiv RSS_{BL}(\boldsymbol{I}_n, \rho(\boldsymbol{X}^\top \boldsymbol{X})^{-1})$ holds if and only if $\boldsymbol{\Omega} = \boldsymbol{I}_n$.
In fact, it follows from \eqref{517} that $(1+\rho)^2 \boldsymbol{\Gamma} = (1+\rho)^2 (\boldsymbol{X}^\top \boldsymbol{X})^{-1} \ \Leftrightarrow \ \boldsymbol{\Gamma} = (\boldsymbol{X}^\top \boldsymbol{X})^{-1}$ under \eqref{516}.
Thus, $\boldsymbol{\Omega}$ is of the form $\boldsymbol{\Omega} = \boldsymbol{X}(\boldsymbol{X}^\top \boldsymbol{X})^{-1}\boldsymbol{X}^\top + \boldsymbol{Z}(\boldsymbol{Z}^\top \boldsymbol{Z})^{-1}\boldsymbol{Z}^\top = \boldsymbol{I}_n$.
\end{rem}

\begin{cor}\label{cor52}
Let $\boldsymbol{\Omega}$ be written as \eqref{Orep}.
For a given $\boldsymbol{K} \in \MS^+(k)$, 
the equality 
\begin{align*}
RSS_{BL}(\boldsymbol{\Omega}, \boldsymbol{K}) \equiv RSS_{BL}(\boldsymbol{I}_n,\boldsymbol{K})
\end{align*}
holds if and only if $\boldsymbol{\Omega}$ is of the form \eqref{516} with $\boldsymbol{\Gamma} \in \MS^+(k)$ satisfying 
\begin{align}\label{newG}
\boldsymbol{K}(\boldsymbol{\Gamma}^{-1}-\boldsymbol{X}^{\top}\boldsymbol{X})\boldsymbol{K}=(\boldsymbol{X}^{\top}\boldsymbol{X})^{-1}-\boldsymbol{\Gamma}.
\end{align}
\end{cor}
\begin{proof}
We directly derive the conditions by using Corollary~\ref{cor23} and Corollary 4.3 of \cite{RefMT25}. 
This completes the proof.
\end{proof}
\begin{rem}
The equality $RSS_{BL}(\boldsymbol{\Omega}, (\boldsymbol{X}^{\top}\boldsymbol{X})^{-1}) \equiv RSS_{BL}(\boldsymbol{I}_n, (\boldsymbol{X}^{\top}\boldsymbol{X})^{-1})$ holds if and only if $\boldsymbol{\Omega} = \boldsymbol{I}_n$.
In fact, it follows from \eqref{newG} that $\Tilde{\boldsymbol{\Gamma}} + \Tilde{\boldsymbol{\Gamma}}^{-1} = 2\boldsymbol{I}_k$, where 
\begin{align*}
\Tilde{\boldsymbol{\Gamma}} = (\boldsymbol{X}^\top \boldsymbol{X})^{1/2} \boldsymbol{\Gamma} (\boldsymbol{X}^\top \boldsymbol{X})^{1/2} \in \MS^+(k).
\end{align*}
From spectral decomposition, there exists an orthogonal matrix $\boldsymbol{P} \in \mathbb{R}^{k \times k}$ such that 
$\Tilde{\boldsymbol{\Gamma}} = \boldsymbol{P}\boldsymbol{\Lambda}\boldsymbol{P}^\top$ with $\boldsymbol{\Lambda} = \diag(\lambda_1, \ldots, \lambda_k)\in \MS^+(k)$.
This implies that $\lambda_i + \lambda_i^{-1} = 2$ for each $i=1, \ldots ,k$, and hence $\lambda_i =1$ for all $i$.
Therefore, we have $\Tilde{\boldsymbol{\Gamma}} = \boldsymbol{I}_k \ \Leftrightarrow \ \boldsymbol{\Gamma} = (\boldsymbol{X}^\top \boldsymbol{X})^{-1}$.
\end{rem}
The following theorem, one of the main results of this paper, provides additional conditions under which the equality \eqref{B2E} guarantees  the trivial condition.
\begin{thm}\label{thm5}
For given $\boldsymbol{K}_1,\boldsymbol{K}_2 \in \MS^N(k)$, the two equalities \eqref{B2E} and \eqref{BRSS1} simultaneously hold if and only if $\boldsymbol{\Omega} = \boldsymbol{I}_n$ and $\boldsymbol{K}_1 = \boldsymbol{K}_2$.
\end{thm}
\begin{proof}
It follows from the equality \eqref{B2E} that $\boldsymbol{y} - \boldsymbol{X}\hat{\boldsymbol{\beta}}_{BL}(\boldsymbol{\Omega},\boldsymbol{K}_1) = \boldsymbol{y} - \boldsymbol{X}\hat{\boldsymbol{\beta}}_{BL}(\boldsymbol{I}_n,\boldsymbol{K}_2)$.
This implies that 
\begin{align*}
&RSS_{BL}(\boldsymbol{\Omega},\boldsymbol{K}_1)=\left(\boldsymbol{y}-\boldsymbol{X}\hat{\boldsymbol{\beta}}_{BL}(\boldsymbol{I}_n,\boldsymbol{K}_2)\right)^{\top}\boldsymbol{\Omega}^{-1}\left(\boldsymbol{y}-\boldsymbol{X}\hat{\boldsymbol{\beta}}_{BL}(\boldsymbol{I}_n,\boldsymbol{K}_2)\right)=\boldsymbol{y}^{\top}\boldsymbol{T}\boldsymbol{\Omega}^{-1}\boldsymbol{T}\boldsymbol{y}
\end{align*}
and $RSS_{BL}(\boldsymbol{I}_n,\boldsymbol{K}_2)=\boldsymbol{y}^{\top}\boldsymbol{T}\boldsymbol{T} \boldsymbol{y}$.
Thus, the equality \eqref{BRSS1} is written as $\boldsymbol{y}^\top (\boldsymbol{T}\boldsymbol{\Omega}^{-1}\boldsymbol{T}-\boldsymbol{T}\boldsymbol{T})\boldsymbol{y} \equiv 0$, which is further equivalent to
\begin{align}\label{OEI}
&\boldsymbol{T}\boldsymbol{\Omega}^{-1}\boldsymbol{T} = \boldsymbol{T}\boldsymbol{T} \notag \\
\Leftrightarrow \quad & \boldsymbol{\Omega} = \boldsymbol{I}_n.
\end{align}
The equality \eqref{B2E} holds if and only if \eqref{B2Enew} holds from Theorem~\ref{thm2}.
Hence, we have
\begin{align*}
&\boldsymbol{\Omega X K}_2 = \boldsymbol{X}\boldsymbol{K}_1\\
\Leftrightarrow \quad &\boldsymbol{X K}_2 = \boldsymbol{X}\boldsymbol{K}_1\\
\Leftrightarrow \quad &\boldsymbol{K}_1 =\boldsymbol{K}_2
\end{align*}
under \eqref{OEI}.
This completes the proof.
\end{proof}

\section{Concluding remarks}\label{sec:7}
In this paper, we investigated statistical properties and equivalence problems of BLEs in the general linear model.
We first examined the linear sufficiency and linear completeness of BLE.
Then we derived necessary and sufficient conditions for the equality of two BLEs and showed that these conditions can simplify estimation procedures when the covariance matrix contains unknown parameters.
Finally, we obtained several results regarding the equivalence of the residual sums of squares.

For future directions, it is interesting to study the case where $\boldsymbol{\Omega}$ is singular or $\boldsymbol{X}$ does not have full column rank. 

\section*{Acknowledgments}
The author would like to express his sincere gratitude to the anonymous reviewers for their valuable comments.
The author is also deeply grateful to his supervisor, Prof. Koji Tsukuda, for his insightful advice.


\begin{thebibliography}{99}	
%
\bibitem[Anderson(1948)]{RefA48}
Anderson, T.~W. (1948).
On the theory of testing serial correlation.
\textit{Skand. Aktuarietidskr.} \textbf{31}, 88--116.
%
\bibitem[Arnold and Stahlecker(2000)]{RefAS00}
Arnold, B.~F., Stahlecker, P. (2000).
Another view of the Kuks--Olman estimator.
\textit{J. Statist. Plann. Inference} \textbf{89}, no.1--2, 169--174.
%
\bibitem[Arnold and Stahlecker(2011)]{RefAS11}
Arnold, B.~F., Stahlecker, P. (2011).
An unexpected property of minimax estimation in the relative squared error approach to linear regression analysis.
\textit{Metrika} \textbf{74}, no.3, 397--407.
%
\bibitem[Baksalary et al.(1989)]{RefB89}
Baksalary, J.~K., Liski, E.~P., Trenkler, G. (1989). Mean square error matrix improvements and admissibility of linear estimators. \textit{J. Statist. Plann. Inference} \textbf{23}, no.3, 313--325.
%
\bibitem[Baksalary and Markiewicz(1988)]{RefBM88}
Baksalary, J.~K., Markiewicz, A. (1988). Admissible linear estimators in the general Gauss-Markov model. \textit{J. Statist. Plann. Inference} \textbf{19}, no.3, 349--359.
%
\bibitem[Bunke(1975)]{RefB75}
Bunke, O. (1975).
Minimax linear, ridge and shrunken estimators for linear parameters.
\textit{Math. Operationsforsch. Statist.} \textbf{6}, no.5, 697--701.
%
\bibitem[Cnaan et al.(1997)]{RefC97}
Cnaan, A., Laird, N. M., Slasor, P. (1997).
Using the general linear mixed model to analyse unbalanced repeated
measures and longitudinal data.
\textit{Stat. Med.} \textbf{16}, 2349--2380.
%
\bibitem[Dow et al.(1982)]{RefD82}
Dow, M. M., Burton, M. L., White, D. R. (1982).
Network autocorrelation: A simulation study of a foundational problem in regression and survey research.
\textit{Social Networks} \textbf{4}, no.2, 169--200.
%
\bibitem[Drygas(1983)]{RefD83}
Drygas, H. (1983).
Sufficiency and completeness in the general Gauss--Markov model.
\textit{Sankhy\={a} Ser.~A} \textbf{45}, no.1, 88--98.
%
\bibitem[Geisser(1970)]{RefG70}
Geisser, S. (1970).
Bayesian analysis of growth curves.
\textit{Sankhy\={a} Ser.~A} \textbf{32}, 53--64.
%
\bibitem[Gotu(2001)]{RefG01}
Gotu, B. (2001). The equality of OLS and GLS estimators in the linear regression model when the disturbances are spatially correlated. \textit{Statist. Papers} \textbf{42}, 253--263.
%
\bibitem[Gro\ss(2003)]{RefG03}
Gro\ss, J. (2003). \textit{Linear Regression}. Springer, Berlin.
%
\bibitem[Gro\ss(2004)]{RefG04}
Gro\ss, J. (2004).
The general Gauss–Markov model with possibly singular dispersion matrix. \textit{Statist. Papers} \textbf{45}, no.3, 311--336.
%
\bibitem[Gro\ss \ and Markiewicz(2004)]{RefGM04}
Gro\ss, J., Markiewicz, A. (2004).
Characterizations of admissible linear estimators in the linear model.
\textit{Linear Algebra Appl.} \textbf{388}, 239--248.
%
\bibitem[Gro\ss \ and Trenkler(1997)]{RefGT97}
Gro\ss, J., Trenkler, G. (1997). When do linear transforms of ordinary least squares and Gauss--Markov estimator coincide?
\textit{Sankhy\={a} Ser.~A} \textbf{59}, no.2, 175--178.
%
\bibitem[Gro\ss \ et al.(2001)]{RefG01b}
Gro\ss, J., Trenkler, G., Werner, H. J. (2001).
The equality of linear transforms of the ordinary least squares estimator and the best linear unbiased estimator.
\textit{Sankhy\={a} Ser.~A} \textbf{63}, no.1, 118--127.
%
\bibitem[Hartigan(1969)]{RefH69}
Hartigan, J. A. (1969). Linear Bayesian methods. \textit{J. Royal Statist. Soc. Ser. B} \textbf{31}, no.3, 446--454.
%
\bibitem[Haslett et al.(2023)]{RefHIMP23}
Haslett, S. J., Isotalo, J., Markiewicz, A., Puntanen, S. (2023).
Further remarks on permissible covariance structures for simultaneous retention of BLUEs in linear models.
\textit{Acta Comment. Univ. Tartu. Math.} \textbf{27}, no.2, 101--112.
%
\bibitem[Haslett et al.(2025)]{RefHIMP25}
Haslett, S.~J., Isotalo, J., Markiewicz, A., Puntanen, S. (2025).
Linear sufficiency and permissible covariance structures for retention of BLUEs in linear models.
\textit{Indian J. Pure Appl. Math.} \textbf{56}, no.3, 901--916.
%
\bibitem[Haslett et al.(2020)]{RefH20}
Haslett, S. J., Liu, X. Q., Markiewicz, A., Puntanen, S. (2020).
Some properties of linear sufficiency and the BLUPs in the linear mixed model.
\textit{Statist. Papers} \textbf{61}, no.1, 385--401.
%
\bibitem[Hauke et al.(2012)]{RefH12}
Hauke, J., Markiewicz, A., Puntanen, S. (2012).
Comparing the BLUEs under two linear models.
\textit{Comm. Statist. Theory Methods} \textbf{41}, no.13--14, 2405--2418.
%
\bibitem[Hoerl and Kennard(1970a)]{RefHK70a}
Hoerl, A. E., Kennard, R. W. (1970a). Ridge regression: biased estimation for nonorthogonal problems. \textit{Technometrics} \textbf{12}, no.1, 55--67.
%
\bibitem[Hoerl and Kennard(1970b)]{RefHK70b}
Hoerl, A. E., Kennard, R. W. (1970b). Ridge regression: applications to nonorthogonal problems. \textit{Technometrics} \textbf{12}, no.1, 69--82.
%
\bibitem[Hoffmann(1996)]{RefH96}
Hoffmann, K. (1996). A subclass of Bayes linear estimators that are minimax. \textit{Acta Appl. Math.} \textbf{43}, no.1, 87--95.
%
\bibitem[Ibarrola and Pérez-Palomares(2003)]{RefIP03}
Ibarrola, P., Pérez-Palomares, A. (2003).
Linear sufficiency and linear admissibility in a continuous time Gauss-Markov model.
\textit{J. Multivariate Anal.} \textbf{87}, no.2, 315--327.
%
\bibitem[Ibarrola and Pérez-Palomares(2004)]{RefIP04}
Ibarrola, P., Pérez-Palomares, A. (2004).
Linear completeness in a continuous time Gauss-Markov model.
\textit{Statist. Probab. Lett.} \textbf{69}, no.2, 143--149.
%
\bibitem[Isotalo and Puntanen(2006)]{RefIP06}
Isotalo, J., Puntanen, S. (2006).
Linear sufficiency and completeness in the partitioned linear model.
\textit{Acta Comment. Univ. Tartu. Math.} \textbf{10}, 53--67.
%
\bibitem[Jiang et al.(2022)]{RefJ22}
Jiang, J., Wang, L., Wang, L. (2022).
Linear approximate Bayes estimator for regression parameter with an inequality constraint.
\textit{Comm. Statist. Theory Methods} \textbf{51}, no.6, 1531--1548.
%
\bibitem[Kala et al.(2017)]{RefK17}
Kala, R., Puntanen, S., Tian, Y. (2017).
Some notes on linear sufficiency.
\textit{Statist. Papers} \textbf{58}, no.1, 1--17.
%
\bibitem[Kariya(1980)]{RefK80}
Kariya, T. (1980).
Note on a condition for equality of sample variances in a linear model.
\textit{J. Amer. Statist. Assoc.} \textbf{75}, no.371, 701--703.
%
\bibitem[Kariya and Kurata(2004)]{RefKK04}
Kariya, T., Kurata, H. (2004). \textit{Generalized Least Squares}. Wiley Series in Probability and Statistics. John Wiley \& Sons, Ltd., Chichester.
%
\bibitem[Kr$\ddot{\mbox{a}}$mer(1980)]{RefK80b}
Kr$\ddot{\mbox{a}}$mer, W. (1980). A note on the equality of ordinary least squares and Gauss--Markov estimates in the general linear model.
\textit{Sankhy\={a} Ser.~A} \textbf{42}, 130--131.
%
\bibitem[Kr$\ddot{\mbox{a}}$mer and Donninger(1987)]{RefKD87}
Kr$\ddot{\mbox{a}}$mer, W., Donninger, C. (1987).
Spatial autocorrelation among errors and the relative efficiency of OLS in the linear regression model.
\textit{J. Amer. Statist. Assoc.} \textbf{82}, no.398, 577--579.
%
\bibitem[Kruskal(1968)]{RefK68}
Kruskal, W. (1968).
When are Gauss--Markov and least squares estimators identical? A coordinate-free approach.
\textit{Ann. Math. Statist.} \textbf{39}, 70--75.
%
\bibitem[Kurata(1998)]{RefK98}
Kurata, H. (1998).
A generalization of Rao's covariance structure with applications to several linear models.
\textit{J. Multivariate Anal.} \textbf{67}, no.2, 297--305.
%
\bibitem[LaMotte(1978)]{RefL78}
LaMotte, L.~R. (1978).
Bayes linear estimators.
\textit{Technometrics} \textbf{20}, no.3, 281--290.
%
\bibitem[Luati and Proietti(2011)]{RefLP11}
Luati, A., Proietti, T. (2011).
On the equivalence of the weighted least squares and the generalised least squares estimators, with applications to kernel smoothing.
\textit{Ann. Inst. Stat. Math.} \textbf{63}, no.4, 851--871.
%
\bibitem[Markiewicz(1996)]{RefM96}
Markiewicz, A. (1996).
Characterization of general ridge estimators.
\textit{Statist. Probab. Lett.} \textbf{27}, no.2, 145--148.
%
\bibitem[Markiewicz et al.(2021)]{RefM21}
Markiewicz, A., Puntanen, S., Styan, G.~P.~H. (2021).
The legend of the equality of OLSE and BLUE: highlighted by C.~R. Rao in 1967.
\textit{Methodology and applications of statistics, Contrib. Stat.}, 51--76, Springer, Cham.
%
\bibitem[Marquardt(1970)]{RefM70}
Marquardt, D. W. (1970). Generalized inverses, ridge regression, biased linear estimation, and
nonlinear estimation. \textit{Technometrics} \textbf{12}, 591--612.
%
\bibitem[Mayer and Willke(1973)]{RefMW73}
Mayer, L.~S., Willke, T.~A. (1973).
On biased estimation in linear models. \textit{Technometrics} \textbf{15}, 497--508.
%
\bibitem[Mitra and Moore(1973)]{RefMM73}
Mitra, S.~K., Moore, B.~J. (1973).
Gauss--Markov estimation with an incorrect dispersion matrix.
\textit{Sankhy\={a} Ser.~A} \textbf{35}, 139--152.
%
\bibitem[Mukasa(2026)]{RefM26}
Mukasa, H. (2026).
Linear Estimators, Bayes. \textit{Wiley StatsRef: Statistics Reference Online}.
Wiley, Chichester.
%
\bibitem[Mukasa and Tsukuda(2025)]{RefMT25}
Mukasa, H., Tsukuda, K. (2025).
Equality between two general ridge estimators and equivalence of their residual sums of squares. \textit{Statist. Papers} \textbf{66}, no.1, 27.
%
\bibitem[Puntanen and Styan(1989)]{RefPS89}
Puntanen, S., Styan, G.~P.~H. (1989).
The equality of the ordinary least squares estimator and the best linear unbiased estimator.
\textit{Amer. Statist.} \textbf{43}, no.3, 153--164.
%
\bibitem[Puntanen et al.(2005)]{RefP05}
Puntanen, S., Styan, G. P. H., Tian, Y. (2005).
Three rank formulas associated with the covariance matrices of the BLUE and the OLSE in the general linear model.
\textit{Econometric Theory} \textbf{21}, no.3, 659--663.
%
\bibitem[Rao(1967)]{RefR67}
Rao, C.~R. (1967).
Least squares theory using an estimated dispersion matrix and its application to measurement of signals.
\textit{Proc. Fifth Berkeley Sympos. Math. Statist. and Probability (Berkeley, Calif., 1965/66), Vol. I: Statistics}, pp.355--372. University of California Press, Berkeley, CA.
%
\bibitem[Rao(1975)]{RefR75}
Rao, C.~R. (1975).
Simultaneous estimation of parameters in different linear models and applications to biometric problems.
\textit{Biometrics} \textbf{31}, 545--553.
%
\bibitem[Rao(1976)]{RefR76}
Rao, C.~R. (1976).
Estimation of parameters in a linear model.
\textit{Ann. Statist.} \textbf{4}, no.6, 1023--1037.
%
\bibitem[Rao et al.(2007)]{RefR07}
Rao, C.~R., Toutenburg, H., Shalabh, Heumann, C. (2007).
\textit{Linear Models and Generalizations:
Least Squares and Alternatives, Third Edition}. Springer Berlin, Heidelberg.
%
\bibitem[Tian and Wiens(2006)]{RefTW06}
Tian, Y., Wiens, D.~P. (2006).
On equality and proportionality of ordinary least squares, weighted least squares and best linear unbiased estimators in the general linear model.
\textit{Statist. Probab. Lett.} \textbf{76}, no.12, 1265--1272.
%
\bibitem[Tsukuda and Kurata(2020)]{RefTK20}
Tsukuda, K., Kurata, H. (2020).
Covariance structure associated with an equality between two general ridge estimators.
\textit{Statist. Papers} \textbf{61}, no.3, 1069--1084.
%
\bibitem[Wang et al.(2024)]{RefW24}
Wang, Q., Wang, L.,  Wang, L. (2024).
Bayesian instrumental variable estimation in linear measurement error models.
\textit{Canad. J. Statist.} \textbf{52}, no.2, 500--531.
%
\bibitem[Zhang and Wang(2021)]{RefZW21}
Zhang, F., Wang, L. (2021).
Linear Bayesian estimators for linear models with constraints.
\textit{J. Stat. Comput. Simul.} \textbf{91}, no.8, 1635--1650.
%
\bibitem[Zyskind(1967)]{RefZ67}
Zyskind, G. (1967).
On canonical forms, non-negative covariance matrices and best and simple least squares linear estimators in linear models.
\textit{Ann. Math. Statist.} \textbf{38}, 1092--1109.
%
\end{thebibliography}
\end{document}